\documentclass[12pt]{amsart}
\usepackage{amsmath}
\usepackage{amsfonts}
\usepackage{amssymb}
\def\N{{\Bbb N}}     
\def\Z{{\Bbb Z}}
\def\Q{{\Bbb Q}}    
\def\R{{\Bbb R}}
\def\C{{\Bbb C}}

\def\Sp{\operatorname{Sp}}
\newcommand{\SL}{{\operatorname{SL}}}
\newcommand{\PSL}{{\operatorname{PSL}}}
\newcommand{\GL}{{\operatorname{GL}}}

\newcommand{\fraca}{\mathfrak a}
\newcommand{\fracf}{\mathfrak f}
\newcommand{\fracg}{\mathfrak g}

\newcommand{\Half}{\mathcal{H}}

\newcommand{\Pn}{\mathcal{P}_n}

\newcommand{\Pnbar}{\bar{\mathcal{P}}_n}
\newcommand{\tr}{{\operatorname{tr}}}

\newcommand{\Aut}{{\operatorname{Aut}}}
\newcommand{\ord}{{\operatorname{ord}}}

\newcommand{\FMpm}{{\mathbb M}_{\pm}}
\newcommand{\FMcusppm}{{\mathbb S}_{\pm}}
\newcommand{\FM}{\mathbb M}
\newcommand{\FMcusp}{\mathbb S}
\newcommand{\XtwoN}{{\mathcal X}(N)}
\newcommand{\XtwosN}{{\bar{\mathcal X}}(N)}

\newcommand{\FJ}{{\operatorname{FJ}}}
\newcommand{\FS}{{\operatorname{FS}}}
\newcommand{\AFS}{{\operatorname{AFS}}}
\newcommand{\inv}{^{-1}}
\def\JkNm{J_{k,Nm}}  

\def\Jkm{J_{k,m}}
\newcommand\Jkmc{J_{k,m}^{\text{\rm cusp}}}
\newcommand{\GammaopmN}{\Gamma^0_{\pm}(N)}
\newcommand\smallmat[4]{\left(\begin{smallmatrix}
{#1}&{#2}\\{#3}&{#4}\end{smallmatrix}\right)}
\newcommand\Spgp{\operatorname{Sp}}
\def\Stab{\operatorname{Stab}}
\newcommand\KN{K(N)}
\newcommand\KNp{K(N)^{+}}
\newcommand{\diag}{\operatorname{diag}}
\newcommand{\IM}[1]{{\rm Im}(#1)}
\newcommand{\tzzw}{\smallmat{\tau}{z}{z}{\omega}}
\newcommand{\nrNm}{\smallmat{n}{r/2}{r/2}{Nm}}
\newcommand{\floor}{\operatorname{floor}}

\newcommand{\ringmono}{\iota}
\newcommand{\HB}{\operatorname{HB}}
\newcommand{\supp}{\operatorname{supp}}
\newcommand{\divv}{\operatorname{div}}
\newcommand{\Cupp}{\operatorname{C_N}}
\newcommand{\cupp}{\operatorname{c_N}}
\newcommand{\ffs}{\C^{\XtwosN}}
\newcommand{\Ptwoone}{\operatorname{P}_{2,1}}
\newcommand{\PtwooneZ}{\Ptwoone(\Z)}
\newcommand{\PtwooneQ}{\Ptwoone(\Q)}

\newcommand{\Sps}{{\mathcal S}}
\newcommand{\Taylor}{T_0}
\newcommand{\Atz}{A}
\newcommand{\Hatz}{H}
\newcommand{\Horo}{N_{\infty}}
\newcommand{\Grit}{\operatorname{Grit}}

\newcommand{\DD}{\mathcal{D}}
\newcommand{\EE}{\mathcal{E}}
\newcommand{\MMM}{\mathcal{M}}
\newcommand{\OO}{\mathcal{O}}
\newcommand{\RR}{\mathcal{R}}
\newcommand{\id}{\operatorname{id}}
\newcommand{\loc}{\operatorname{loc}}
\newcommand{\germp}{\operatorname{germ}_p}
\newcommand\jj{\operatorname{j}}
\newcommand\WW{\mathcal W}
\newcommand\Denom{\operatorname{Denom}}
\newcommand\Pj{\mathbb P}
\newcommand{\Ponezero}{\mathcal{P}_{10}}
\newcommand\Wpoly{W}

\usepackage{color}

\newtheorem{theorem}{Theorem}[section]
\newtheorem{lemma}[theorem]{Lemma}
\newtheorem{prop}[theorem]{Proposition}
\newtheorem{cor}[theorem]{Corollary}

\newtheorem{defn}[theorem]{Definition}

\allowdisplaybreaks

\begin{document}
\date{\today}
\title[Formal Series]{Formal series of Jacobi forms}
\author[Aoki, Ibukiyama, Poor]{Hiroki Aoki, Tomoyoshi Ibukiyama, Cris Poor}

\address{Department of Mathematics, Faculty of Science and Technology, 
Tokyo University of Science, Noda, Chiba, 278-8510 Japan}
%\email{aoki_hiroki_math@nifty.com}
\email{aoki$\_\hskip0.015in$hiroki$\_\hskip0.015in$math@nifty.com}

\address{Department of Mathematics, Graduate School of Mathematics, Osaka University,
Machikaneyama 1-1, Toyonaka, Osaka, 560-0043 Japan}
\email{ibukiyam@math.sci.osaka-u.ac.jp}

%\author[C.~Poor]{Cris Poor}
\address{Department of Mathematics, Fordham University, Bronx, NY 10458 USA}
\email{poor@fordham.edu}

\begin{abstract}
We prove for general paramodular level 
that formal series of scalar Jacobi forms with an involution condition necessarily converge 
and are therefore the Fourier-Jacobi expansions at the standard 1-cusp 
of paramodular Fricke eigenforms.  
\end{abstract}

%\subjclass[2020]{Primary: 11F46; secondary: 11F30,11F50}

\maketitle

\tableofcontents

%%%%%%%%%%%%%%%%%%%%%%%%%%%
%%%%%%%%%%%%%%%%%%%%%%%%%%%
\section{Introduction.}
\label{secintro}
%%%%%%%%%%%%%%%%%%%%%%%%%%%
%%%%%%%%%%%%%%%%%%%%%%%%%%%

A parmodular form is a Siegel modular form of degree two 
for the following discrete group
$$
K(N)= \Spgp(4,\Q) \cap
\begin{pmatrix} 
\Z  &  N \Z  &  \Z  &  \Z  \\
\Z  &   \Z  &  \Z  &  \frac1{N}\Z  \\
\Z  &  N \Z  &  \Z  &  \Z  \\
N\Z  &  N \Z  &  N\Z  &  \Z  
\end{pmatrix}
=\Stab_{\Spgp(4,\Q)}
\begin{pmatrix} 
\Z  \\  \Z  \\  \Z  \\  N\Z  
\end{pmatrix}.  
$$
Paramodular forms $M_k\left( \KN \right)$ 
are a natural generalization of elliptic modular forms $M_k\left( \Gamma_0(N) \right)$ and 
are interesting in many ways.  
Roberts and Schmidt~\cite{robertsschmidt06,robertsschmidt07}  
gave a sophisticated theory of local and global paramodular newforms.  
Paramodular newforms have applications to modularity; weight two 
to abelian surfaces~\cite{brkr14,brkr18}, and weight three to nonrigid Calabi-Yau threefolds~\cite{gvs} 
of Hodge type~$(1,1,1,1)$.  
Conjectures of Ibukiyama~\cite{i85,i07b}, and of 
Ibukiyama and Kitayama~\cite{ik17},  
connecting paramodular forms to algebraic modular forms 
on a compact twist of $\operatorname{GSp}(4)$, 
motivated mainly by independent calculations of dimension formulae, 
have recently been proven~\cite{MR4329279,MR4751186}. 
In~\cite{MR4730249}, this connection was generalized to a correspondence 
with Fricke eigenspaces.  
Utilizing results of~\cite{MR4730249},  
dimension formulae for Fricke plus and minus spaces 
of paramodular forms for prime level were computed in~\cite{ibuk24}.  
For weights~$k \ge 3$, these recent proofs allow paramodular Hecke eigensystems to also be computed using
orthogonal modular forms~\cite{MR4730249,MR4732690}.  
Finally, Gritsenko lifts and Borcherds products provide 
concrete examples of paramodular forms~\cite{gn98}.   
\smallskip 

Paramodular forms $f \in M_k\left( \KN \right)$ 
have Fourier expansions in three variables but the natural 
generalization of the Fourier series of an elliptic modular form is perhaps the 
Fourier-Jacobi expansion of a paramodular form, 
which, using the components $\Omega=\tzzw \in \Half_2$, 
recollects the Fourier series in powers of~$e(\omega)=e^{2\pi i \omega}$,    
$$
f\tzzw = \sum_{m=0}^{\infty} 
\phi_m(\tau,z) e\left( Nm\omega \right).  
$$
Each coefficient is a Jacobi form, $\phi_m \in \JkNm$,  
and has its own Fourier expansion 
$\phi_m(\tau,z) = \sum_{n \ge 0, r \in \Z} c(n,r;\phi_m) e(n\tau+r z)$.    
The Fourier-Jacobi expansion thus defines a map to formal series of Jacobi forms 
$\FJ: M_k\left( \KN \right) \to \FM(k,N) = \prod_{m=0}^{\infty} \JkNm$.  
This map is not surjective because we cannot freely select a sequence of Jacobi forms and 
obtain the convergent Fourier-Jacobi expansion of a paramodular form.  
One source of consistency conditions among the Fourier-Jacobi coefficients~$(\phi_m)$ 
arises from a normalizing involution of~$\KN$, 
the paramodular Fricke involution 
$\mu_N = \smallmat{{{}^tF_N\inv}}{0}{0}{F_N}$, where
$F_N= {\tiny\frac{1}{\sqrt{N}}} \smallmat01{-N}0$ is the Fricke involution on~$\Gamma_0(N)$.  
This involution splits $M_k\left( \KN \right)$ into plus and minus forms, 
$M_k\left( \KN \right)=M_k\left( \KN \right)^{+} \oplus M_k\left( \KN \right)^{-}$.  
The block diagonal form of~$\mu_N$ gives a simple action on the Fourier series and 
consequently gives the following 
{\em involution condition\/} on the Fourier-Jacobi coefficients of any 
$f \in M_k\left( \KN \right)^{\epsilon}$, $\epsilon \in \{ \pm 1\}$:    
\begin{align}
\label{eqinvcondintro}
%\notag
&\text{For all semidefinite  $\smallmat{n}{r/2}{r/2}{Nm}$ with $n,r,m\in \Z$, }  \\  \notag
&\hskip0.7in c(n,r; \phi_{m})= (-1)^k \epsilon\, c(m,r; \phi_{n}).  
\end{align}
Let $\FM(k,N,\epsilon)$ be the subspace of formal series of Jacobi forms 
satifying the involution condition:  
$$
\FM(k,N,\epsilon)= 
\{ \fracf  \in \FM(k,N): 
\text{\rm $\fracf$ satisfies condition~{(\ref{eqinvcondintro})}}\}.  
$$
The Fourier-Jacobi expansion gives $\FJ: M_k\left( \KN \right)^{\epsilon} \to \FM(k,N,\epsilon)$.  
It seems a bit audacious to hope that the involution condition alone 
forces the convergence of a formal series to a paramodular Fricke eigenform, but 
theoretical results for low level~$N$ and computed examples suggest that this is true, 
and we prove it here.

\begin{theorem}[Main Theorem] 
\label{thmain3intro}
Let $N\in\N$, $k \in \N_0$, and $\epsilon \in \{ \pm 1\}$. 
The map 
$\FJ: M_k\left( \KN \right)^{\epsilon} \to \FM(k,N,\epsilon)$ 
is an isomorphism.     
\end{theorem}

The case~$k=0$ is not difficult, and neither is the case~$k=1$, when 
both the domain and codomain are trivial 
by a result of Skoruppa~\cite{MR0806354}, 
but~$k=2$ already has important examples.   
The case~$N=1$ of Theorem~\ref{thmain3intro} was proven 
by the first author~\cite{MR1794522},  
the cases~$2\le N \le 4$ 
by Yuen and the second and third authors~\cite{ipy13}.   
These results raised the question, made explicit in~\cite{ipy13}, of whether 
the involution condition alone implies the convergence of a formal series of Jacobi forms.   
Theorem~\ref{thmain3intro} resolves this question in the affirmative for the first time.  
\smallskip 

Bruinier~\cite{MR3283640}, Raum~\cite{MR3433884}, and Bruinier and Raum~\cite{MR3406827,bruinierraum24}  
have approached the theory of formal series of Jacobi forms in a more general setting but with a different set 
of hypotheses.  
Bruinier~\cite{MR3283640} and Raum~\cite{MR3433884} independently extended Aoki's result for 
$\Sp(4,\Z)$ to vector valued formal series of Jacobi forms, and thereby resolved a conjecture of 
Kudla~\cite{MR1427845,MR1988501,MR2083214} concerning generating series of special cycles for degree~$n=2$.  
In~\cite{MR3406827},  Bruinier and Raum proved that formal series of Jacobi forms for 
$\Sp(2n,\Z)$ converge under the assumption of all $\GL(n,\Z)$ symmetries,  thereby proving  
Kudla's modularity conjecture for general~$n$.  
Pollack \cite{pollack24} has also given a proof of automatic convergence for cuspidal 
automorphic forms on $\Sp(2n)$ that takes Jacobi and Levi symmetries as hypotheses.  
Bruinier and Raum~\cite{bruinierraum24} have recently developed a theory of formal series 
for subgroups commensurable with $\Sp(2n,\Z)$; they prove that 
compatible formal series of Jacobi forms at every $1$-cusp and Levi symmetries 
imply convergence.  
They interpret formal series as sections of a line bundle over a formal complex space.  
Our Main Theorem~\ref{thmain3intro} for paramodular groups 
in degree two and levels $N>3$ 
is not a consequence of any of 
these results.  We only use formal series of Jacobi forms at the standard $1$-cusp and only 
assume symmetry under a single involution.  
\vskip0.1in

Our main result has applications to computations. 
Jacobi restriction~\cite{ipy13,bpy16,MR3713095} is a method that attempts to rigorously 
compute the space $\FJ\left( S_k\left( \KN \right)^{\epsilon}  \right)$ by 
imposing necessary linear relations on the space of Jacobi forms $\prod_{m=1}^{d} \JkNm^{\rm cusp}$. 
Jacobi restriction has provided rigorous upper bounds 
for $\dim S_k\left( \KN \right)^{\epsilon}$, even though the authors could not 
guarantee in advance that the method would work.  
The following corollary proves that any schema for computing spaces of 
paramodular forms that spans spaces of Jacobi forms and imposes the 
involution condition is in principle sound.  

\begin{cor}
\label{corintro}
For $d \in \N$, define the  $\C$-vector space $\FM(k,N,\epsilon)[d]$ as   
$\{ (\phi_m) \in \prod_{m=0}^{d} \JkNm:\text{$(\phi_m)$ satisfies~{\rm (\ref{eqinvcondintro})} for all $n,m \le d$ } \}$.  
The sequence $\dim_{\C} \FM(k,N,\epsilon)[d]$ is monotonically decreasing for $d \ge \frac16 Nk$, 
and 
we have $\lim_{d \to +\infty} \dim_{\C} \FM(k,N,\epsilon)[d] = \dim_{\C} M_k\left( \KN \right)^{\epsilon}$.  
\smallskip

In particular,  we have 
$\dim_{\C} M_k\left( \KN \right)^{\epsilon} \le \dim_{\C} \FM(k,N,\epsilon)[d]$ for $d > \frac16 Nk$, 
and we have equality 
for sufficiently large~$d$.  
\end{cor}

The first author was supported by JSPS KAKENHI Grant Numbers  JP19K03429 and JP23K03039.  
The second author was supported by the following grants: JSPS KAKENHI Grant Number JP19K03424, JP20H00115 and JP23K03031.  
The authors thank the American Institute of Mathematics for its critical support of this research.  

%%%   Hiroki Aoki
%  I was supported by JSPS KAKENHI Grant Number  JP19K03429 and JP23K03039.

%%%   Tomoyoshi Ibukiyama 
%  I was supported by the following grants: JSPS KAKENHI Grant Number JP19K03424, JP20H00115, JP23K03031.

%%%%%%%%%%%%%%%%%%%%%%%%%%%
%%%%%%%%%%%%%%%%%%%%%%%%%%%
\section{Notation.}
\label{secnotation}
%%%%%%%%%%%%%%%%%%%%%%%%%%%
%%%%%%%%%%%%%%%%%%%%%%%%%%%

We denote the natural numbers by $\N=\{1,2,3,\ldots\}$ and the whole numbers 
by $\N_0=\{0,1,2,\ldots\}$.  
Throughout this article $N\in\N$ denotes a level, $k \in \N_0$ a weight, and 
$\epsilon \in \{ \pm 1 \}$ a sign.  
We write $\sigma'$ for the transpose of a matrix~$\sigma$, and 
$\sigma^*$ for the transpose inverse.  
Let $t[\sigma]=\sigma' t \sigma$ for compatibly sized matrices, and 
$\langle a,b \rangle=\tr(ab)$ for $a,b \in M_{n \times n}^{\text{\rm sym}}(\C)$.  
For $z \in \C$, set $e(z)=e^{2\pi i z}$.  
\smallskip

For the theory of Jacobi forms see~\cite{ez85}.  
The upper half plane is $\Half_1=\{ \tau \in \C: \IM \tau >0 \}$.  
For a Jacobi form $\phi \in \Jkm$ 
of weight~$k$ and index~$m \in \N_0$ 
with $\phi: \Half_1 \times \C \to \C$ 
we write the Fourier expansion as 
$\phi(\tau,z) = \sum_{n \in \N_0,\, r \in \Z} c(n,r; \phi) q^n \zeta^r$
with $q=e(\tau)$ and $\zeta=e(z)$.  
The order of a nonzero $\phi$ is 
$\ord \phi = \min\{n \in \N_0: \exists r \in \Z: c(n,r;\phi) \ne 0  \}$.  
For $\nu \in \N_0$, we set $\Jkm(\nu)=\{\phi \in \Jkm: \ord \phi \ge \nu \}$.  
\smallskip

For the theory of Siegel modular forms we refer to~\cite{freitag83}.  
For a ring $R \subseteq \R$ define the positive definite cone 
with entries in~$R$ by 
$\Pn(R)=\{ s \in M_{n\times n}^{\text{\rm sym}}(R): s>0\}$, 
and the positive semidefinite cone by 
$\Pnbar(R)=\{ s \in M_{n\times n}^{\text{\rm sym}}(R): s \ge 0\}$.  
The Siegel upper half space is 
$\Half_n=\{\Omega \in M_{n\times n}^{\text{\rm sym}}(\C): \IM \Omega \in \Pn(\R) \}$.  
An element  $\sigma = \smallmat{a}{b}{c}{d}$ in the real symplectic group 
$\Sp(2n,\R)$ acts on the symmetric space $\Half_n$ by 
$\sigma\langle \Omega \rangle=(a\Omega+b)(c\Omega+d)\inv$, and 
the Siegel factor of automorphy $j:\Sp(2n,\R) \times \Half_n \to \C^{\times}$ is 
$j(\sigma,\Omega)=\det(c\Omega+d)$. 
For a function $f: \Half_n \to \C$ the slash $k$-action 
$(f \vert_k \sigma)(\Omega)=j(\sigma,\Omega)^{-k} f\left( \sigma\langle \Omega \rangle \right)$ 
defines another function $f \vert_k \sigma: \Half_n \to \C$.  
For a discrete group~$\Gamma \subseteq \Sp(2n,\R)$ commensurable with $\Sp(2n,\Z)$, 
write $M_k\left( \Gamma \right)$ for the $\C$-vector space of Siegel modular forms 
of weight~$k$, and $S_k\left( \Gamma \right)$ for the subspace of cusp forms.  
\smallskip

For the paramodular group~$\KN$ in degree~$n=2$, write the Fourier expansion for 
$f \in M_k\left( \KN \right)$ as 
$f(\Omega) = \sum_{t \in \XtwosN} a(t;f) e\left( \langle \Omega, t \rangle \right)$, 
where  $\XtwosN = \{ \smallmat{n}{r/2}{r/2}{Nm} \in {\bar{\mathcal P}}_2(\Q): n,r,m \in \Z  \}$. 
For the Fourier-Jacobi expansion of~$f$, we write $\Omega=\tzzw$ and collect the 
Fourier expansion in~$e(\omega)$ to obtain 
$f\tzzw = \sum_{m=0}^{\infty} \phi_m(\tau,z) e(Nm\omega)$ 
with Fourier-Jacobi coeffcients $\phi_m \in \JkNm$.  
The paramodular group has a normalizing involution, the paramodular Fricke involution, 
given by $\mu_N = \smallmat{{{}^tF_N\inv}}{0}{0}{F_N} \in \Sp(4,\R) $, with 
 $F_N= \frac{1}{\sqrt{N}} \smallmat01{-N}0 \in \SL(2, \R)$.   
Define the Fricke paramodular group $\KN^{+}=\langle K(N), \mu_N \rangle$.  
For $\epsilon \in \{ \pm 1 \}$, let  
$M_k\left( \KN \right)^{\epsilon}=\{ f \in M_k\left( \KN \right): f \vert_k \mu_N= \epsilon f  \}$ 
be the Fricke eigenspaces.  
There are  graded rings,   
$M\left( \KN \right)= \oplus_{k\in \N_0} M_k\left( \KN \right)$, 
$M_{\pm}\left( \KN \right)= \oplus_{k\in \N_0, \epsilon \in \{ \pm 1\}} M_k\left( \KN \right)^{\epsilon}$, 
as well as the graded ring of Fricke plus forms 
$M\left(\KN^{+} \right)= \oplus_{k\in \N_0} M_k\left( \KN^{+}\right)$.  
\smallskip

\noindent
$\bullet$ $\Horo(\eta)=\{ \tau \in \Half_1: \IM \tau > \eta\}$ for $\eta \ge 0$.    
\smallskip 

\noindent 
$\bullet$ 
$
\Ptwoone(\Z)= \Spgp(4,\Z) \cap
\begin{pmatrix} 
\Z  &  0  &  \Z  &  \Z  \\
\Z  &   \Z  &  \Z  &  \Z  \\
\Z  & 0  &  \Z  &  \Z  \\
0  &  0  &  0  &  \Z  
\end{pmatrix}.   
$

%% $\XtwoN = \{ \smallmat{n}{r/2}{r/2}{Nm} \in {\mathcal P}_2(\Q): n,r,m \in \Z  \}$ 

%%%%%%%%%%%%%%%%%%%%%%%%%%%
%%%%%%%%%%%%%%%%%%%%%%%%%%%
\section{Formal Series.}
\label{secformal}
%%%%%%%%%%%%%%%%%%%%%%%%%%%
%%%%%%%%%%%%%%%%%%%%%%%%%%%

A paramodular form $f \in M_k\left( \KN \right)$ has a Fourier series 
$
f(\Omega)=\sum_{t \in \XtwosN} a(t;f) e\left( \langle \Omega,t \rangle \right)
$ 
that converges absolutely and uniformly on compact subsets of~$\Half_2$.  
The absolute convergence allows rearrangement 
into a Fourier-Jacobi series 
$
f(\Omega)=\sum_{m=0}^{\infty} \phi_m(\tau,z) e\left( Nm\omega\right)
$ 
where each Fourier-Jacobi coefficient is a Jacobi form $\phi_m\in \JkNm$.  
This follows from the fact that the Fourier-Jacobi expansion is term-by-term invariant 
under $\PtwooneZ \subseteq \KN$, see~\cite{MR1336601}.    
We define the formal Fourier-Jacobi series $\FJ(f)$ of a paramodular form~$f$ by 
\begin{align*}
\FJ: M_k\left( \KN \right)   &\to \FM(k,N) =\prod_{m=0}^{\infty} \JkNm \\
f &\mapsto \sum_{m=0}^{\infty} \phi_m \xi^{Nm},   
\end{align*}
where~$\xi$ is a place-holding variable.  
The Fourier coefficients of a paramodular form also satisfy symmetries determined by 
the following group:   
%% the equivalence class group of~$\KN$:   
\begin{align*}
%% u_o\left( \KN \right)&= 
\{ \sigma \in \GL(2,\R): \smallmat{\sigma}00{\sigma^*}\in \KN\}   
&= \langle \Gamma^0(N), \diag(1,-1) \rangle 
=\GammaopmN.  
\end{align*}
We call equation~{(\ref{eqgammaoNsym})} below the $\Gamma^0(N)$-symmetries,  
\begin{equation}
\label{eqgammaoNsym}
\forall \sigma \in   \GammaopmN, 
\forall t \in \XtwosN, \, 
a\left( t[\sigma];f\right) = \det(\sigma)^k a\left( t;f\right) .  
\end{equation}
The Fourier coefficients of~$f$ are related to the 
Fourier coefficients of the Jacobi forms~$\phi_m$ 
in the Fourier-Jacobi expansion of~$f$ by 
$$
a\left( \nrNm ;f\right) = c(n,r;\phi_m).  
$$ 
Accordingly, we define a subspace of~$\FM(k,N)$ that satisfies corresponding  
symmetries.  Let $\Gamma \subseteq \GammaopmN$ be a subgroup.   
\begin{align}
\FM(k,N;\Gamma)&= 
\{ \fracf=\sum_{m=0}^{\infty} \phi_m \xi^{Nm} \in \FM(k,N): 
\text{\rm $\fracf$ satisfies equation~{(\ref{eqgammaoNsymformal})}}\},  \notag \\   
\label{eqgammaoNsymformal}
& 
\forall \sigma \in   \Gamma, 
\forall t_1=\smallmat{n_1}{r_1/2}{r_1/2}{Nm_1}, t_2= \smallmat{n_2}{r_2/2}{r_2/2}{Nm_2}\in \XtwosN, \\    
&\text{if } t_1[\sigma]=t_2 \text{ then } c(n_2,r_2; \phi_{m_2})= 
\det(\sigma)^k c(n_1,r_1; \phi_{m_1}).  \notag
\end{align}
Thus we have $\FJ: M_k\left( \KN \right)   \to \FM(k,N;\Gamma)$.  
\smallskip

We write $\phi_m=\phi_m(\fracf)$ if we need to indicate the dependence of this Jacobi form on the 
formal series~$\fracf$.  
It is often helpful to reformulate the definition of $\FM(k,N;\Gamma)$ by using the notation 
$a(t;\fracf)=c(n,r;\phi_m(\fracf))$ for $t=\nrNm \in \XtwosN$.  
In this way equation~{(\ref{eqgammaoNsymformal})} shortens to 
\begin{equation}
\label{eqgammaoNsymuseful}
\forall \sigma \in   \Gamma, 
\forall t \in \XtwosN, \, 
a\left( t[\sigma];\fracf\right) = \det(\sigma)^k a\left( t;\fracf\right) .  
\end{equation}
We may sometimes avoid tracking boundary conditions in summations by setting $a\left( t;\fracf\right)=0$ for 
$t \in M_{2\times 2}^{\rm sym}(\Q) \setminus \XtwosN$.  
\smallskip

\begin{lemma}
\label{lemringone} %%%  
 Let $\Gamma \subseteq \GammaopmN$ be a subgroup.  
The Cauchy product gives $\FM(N)=\oplus_{k=0}^{\infty} \FM(k,N)$ the structure of a graded ring,  
and $\FM(N;\Gamma)=\oplus_{k=0}^{\infty} \FM(k,N;\Gamma)$ a graded subring.  In particular, 
$\FM(k_1,N) \FM(k_2,N) \subseteq \FM(k_1+k_2,N)$, and 
$\FM(k_1,N;\Gamma) \FM(k_2,N;\Gamma) \subseteq \FM(k_1+k_2,N;\Gamma)$.  
The sets $\FMcusp(N)=\oplus_{k=0}^{\infty} \prod_{m=1}^{\infty} \JkNm^{\rm cusp}$, and 
$\FMcusp(N;\Gamma)=\FM(N;\Gamma) \cap \FMcusp(N)$ are graded ideals in $\FM(N)$, and $\FM(N;\Gamma)$, respectively.  
The Fourier-Jacobi map $\FJ: M\left( \KN \right)   \to \FM(N;\Gamma)$ is an injective ring homomorphism of graded rings 
that sends $S\left( \KN \right)$ to $\FMcusp(N;\Gamma)$, and the inverse image of~$\FMcusp(N;\Gamma)$ 
is $S\left( \KN \right)$.  
\end{lemma}
\begin{proof}
When the formal series $\fracf_j \in \FM(k_j,N)$, for $j=1,2$, are written as 
$\fracf_j = \sum_{m=0}^{\infty} \phi_m(\fracf_j) \xi^{Nm}$ 
then by the definition of the Cauchy product we have 
$\fracf_1 \fracf_2=  \sum_{m=0}^{\infty} \phi_m(\fracf_1 \fracf_2) \xi^{Nm}$ where 
$$
\phi_m(\fracf_1 \fracf_2) =\sum_{m_1,m_2\in \N_0:\, m_1+m_2=m} \phi_{m_1}(\fracf_1)\phi_{m_2}(\fracf_2).  
$$
Thus $\FM(k_1,N) \FM(k_2,N) \subseteq \FM(k_1+k_2,N)$ follows from the grading 
$J_{k_1,N m_1} J_{k_2,N m_2} \subseteq J_{k_1+k_2,N (m_1+m_2)}$ on Jacobi forms, which shows that    
$\FM(N)$ is a graded ring.  To derive a similar result for~$\FM(N;\Gamma)$ 
we reformulate the Cauchy product on~$\FM(N)$ as 
\begin{equation}
\label{eqcauchyprod}
\forall t \in \XtwosN, \, 
a(t;\fracf_1 \fracf_2) =\sum_{t_1,t_2\in \XtwosN:\, t_1+t_2=t}
a(t_1;\fracf_1)a(t_2;\fracf_2).  
\end{equation}
To see that~{(\ref{eqcauchyprod})} follows from the Cauchy product, 
note the following four equalities.  
\begin{align*}
&a(t;\fracf_1 \fracf_2) = c(n,r;\phi_m\left(\fracf_1 \fracf_2)\right) =
\sum_{m_1+m_2=m}c\left(n,r; \phi_{m_1}(\fracf_1)\phi_{m_2}(\fracf_2)\right) \\
= 
&\sum_{\substack{m_1, m_2 \in \N_0:  \\ m_1+m_2=m }} \,
\sum_{\substack{n_1, n_2 \in \N_0:  \\ n_1+n_2=n }}\, 
\sum_{\substack{r_1, r_2 \in \Z:  \\ r_1+r_2=r }} 
c\left(n_1,r_1;\phi_{m_1}(\fracf_1)\right) c\left(n_2,r_2;\phi_{m_2}(\fracf_2)\right) \\
=&\sum_{t_1+t_2=t}
a(t_1;\fracf_1)a(t_2;\fracf_2). 
\end{align*}
Conversely, equation~{(\ref{eqcauchyprod})} implies the second equality, 
which is equivalent to the Cauchy product.  
If we assume $\fracf_j \in \FM(k_j,N;\Gamma)$ for $j=1,2$, then, for any $\sigma \in \Gamma$, 
$a(t[\sigma];\fracf_1 \fracf_2) =\sum_{t_1+t_2=t[\sigma]}
a(t_1;\fracf_1)a(t_2;\fracf_2)$.  
We use 
$\{ (t_1,t_2) \in \XtwosN^2: t_1{+}t_2{=}t[\sigma] \}=\{ (s_1[\sigma],s_2[\sigma]) \in \XtwosN^2: s_1{+}s_2{=}t \}$ 
to change the indices of summation, 
\begin{align*}
a(t[\sigma];\fracf_1 \fracf_2) 
&=\sum_{s_1+s_2=t}
 a(s_1[\sigma];\fracf_1)   a(s_2[\sigma];\fracf_2) \\
&=\sum_{s_1+s_2=t}
\det(\sigma)^{k_1} a(s_1;\fracf_1)  \det(\sigma)^{k_2} a(s_2;\fracf_2) \\
&=\det(\sigma)^{k_1+k_2} a(t;\fracf_1 \fracf_2).   
\end{align*} 
Thus $\FM(k_1,N;\Gamma) \FM(k_2,N;\Gamma) \subseteq \FM(k_1+k_2,N;\Gamma)$ and 
$\FM(N;\Gamma)$ is a graded ring, 
noting that a sum of paramodular forms of distinct weights is zero as a holomorphic function 
if and only if each summand is zero as a holomorphic function.    
The sets $\FMcusp(N)$, and 
$\FMcusp(N;\Gamma)$ are graded ideals simply because 
$J_{k_1,N m_1}^{\rm cusp} J_{k_2,N m_2} \subseteq J_{k_1+k_2,N (m_1+m_2)}^{\text{\rm cusp}}$.  
The convergence of the Fourier-Jacobi expansion shows that~$\FJ$ is injective on each graded piece and 
hence on~$M\left( \KN \right)$.  The absolute convergence of the Fourier-Jacobi expansion 
for $f_j \in M_{k_j}\left( \KN \right)$ proves that $\FJ: M\left( \KN \right)   \to \FM(N;\Gamma)$ is a homomorphism 
because for $\Omega=\tzzw \in \Half_2$, 
\begin{align*}
&f_1(\Omega) f_2(\Omega) = \\
&\sum_{m=0}^{\infty} 
\left( \sum_{m_1+m_2=m} \phi_{m_1}(\FJ(f_1))(\tau,z)\, \phi_{m_2}(\FJ(f_2))(\tau,z)  \right) 
e(Nm\omega),
\end{align*}
and therefore 
\begin{align*}
\FJ(f_1 f_2) &= 
\sum_{m=0}^{\infty} 
\left( \sum_{m_1+m_2=m} \phi_{m_1}(\FJ(f_1))\, \phi_{m_2}(\FJ(f_2))  \right) 
\xi^{Nm}  \\
&=\FJ(f_1) \FJ(f_2).  
\end{align*}
For a cusp form $ f \in S\left( \KN \right) $,   the Siegel~$\Phi$ map gives the containment  
$\supp(f)=  
\{ t \in \XtwosN: a(t;f) \ne 0 \} \subseteq \XtwoN$.  By examining the support of~$\phi_m$,  
\begin{align*}
\supp(\phi_m)&= 
\{ (n,r) \in \N_0 \times \Z: c(n,r;\phi_m) \ne0\}  \\ 
&= 
\{ (n,r) \in \N_0 \times \Z: a\left( \nrNm ;f\right) \ne 0\} \\
&\subseteq 
\{ (n,r) \in \N \times \Z: 4Nmn - r^2 >0\}, 
\end{align*}
we see that $\phi_m \in J_{k,N m}^{\text{\rm cusp}}$, and that $\FJ(f) \in \FMcusp(N)$, 
from which $\FJ(f) \in \FMcusp(N;\Gamma)$ follows.  
The assertion that 
$\{ f \in M_k \left(\KN \right): \FJ(f) \in \FMcusp(N;\Gamma)\} = S_k \left(\KN \right)$ 
relies on special properties of the paramodular group.  
For general subgroups commensurable with $\Sp(4,\Z)$ it is not true that the   
Fourier-Jacobi expansion of a Siegel modular form with coefficients that are all cusp forms 
necessarily comes from a  
Siegel modular cusp form; however, the inference is valid for paramodular forms 
because Reefschl{\"a}ger's 
double coset decomposition~\cite{reefschlager73}
of $\KN \backslash \Sp(4,\Q)/ \PtwooneQ$ has representatives of the form 
$\smallmat{u}{0}{0}{u^{*}}$ with $u \in \GL(2,\Q)$; see also Corollary~{2.5} of 
Gritsenko~\cite{MR1336601}.  
\end{proof}
\smallskip

The involution~$\mu_N=\smallmat{F_N^*}00{F_N}$ decomposes paramodular forms into plus and minus forms 
$M_k \left(\KN \right)= M_k \left(\KN \right)^{+} \oplus M_k \left(\KN \right)^{-}$ 
where 
$M_k \left(\KN \right)^{\epsilon} =\{ f \in M_k \left(\KN \right): f \vert_k \mu_N =\epsilon f \}$.  
The action of $F_N^{*}$ on $\XtwosN$ is 
$F_N\inv t F_N^{*}  = \smallmat{m}{-r/2}{-r/2}{Nn}$ if $t= \nrNm \in \XtwosN$.  
Therefore the Fourier coefficients of a Fricke eigenform $f \in M_k \left(\KN \right)^{\epsilon}$ 
satisfy the {\em involution condition\/} 
\begin{equation}
\label{eqinvcond}
\forall t \in \XtwosN, \,
a\left( F_N\inv t F_N^{*};f \right) =  \epsilon\, a\left( t;f \right) .  
\end{equation}
Accordingly, we define a subspace of formal series satisfying the corresponding  
involution condition 
\begin{align}
\FM(k,N,\epsilon)&= 
\{ \fracf=\sum_{m=0}^{\infty} \phi_m \xi^{Nm} \in \FM(k,N): 
\text{\rm $\fracf$ satisfies equation~{(\ref{eqinvcondformal})}}\},  \notag \\   
\label{eqinvcondformal}
& 
\forall \smallmat{n}{r/2}{r/2}{Nm}\in \XtwosN, \,    
 c(n,r; \phi_{m})= (-1)^k \epsilon\,
c(m,r; \phi_{n}).  
\end{align}
The map $\FJ: M_k\left( \KN \right)^{\epsilon}   \to \FM(k,N,\epsilon)$ injects.    
It is helpful to rewrite equation~{(\ref{eqinvcondformal})} as 
$ \forall t\in \XtwosN, \, a\left( t[F_N^{*}];\fracf \right) = (-1)^k \epsilon\, a\left( t;\fracf \right)$.  
If we write $\FM(k,N,\epsilon;\Gamma)=\FM(k,N;\Gamma) \cap \FM(k,N,\epsilon)$ we also have 
$\FJ: M_k\left( \KN \right)^{\epsilon}   \to \FM(k,N,\epsilon;\Gamma)$.  
We set $\FMcusp(k,N,\epsilon)=\FM(k,N,\epsilon) \cap \FMcusp(k,N)$, and 
the more frequently used 
$\FMcusp(k,N,\epsilon;\Gamma)=\FM(k,N,\epsilon;\Gamma) \cap \FMcusp(k,N)$.  
\smallskip
 
\begin{lemma} 
\label{lemringtwo}
 Let $\Gamma \subseteq \GammaopmN$ be a subgroup.    
The Cauchy product gives $\FMpm(N)=\oplus_{k,\epsilon} \FM(k,N,\epsilon)$, 
  $\FMpm(N;\Gamma)=\oplus_{k,\epsilon} \FM(k,N, \epsilon;\Gamma)$, and 
$\FM(N,+)=\oplus_{k=0}^{\infty} \FM(k,N,+)$
the structure of graded rings.    
In particular, we have 
%%  the containments 
$\FM(k_1,N, \epsilon_1) \FM(k_2,N, \epsilon_2) \subseteq \FM(k_1+k_2,N, \epsilon_1 \epsilon_2)$, and 
$\FM(k_1,N, \epsilon_1;\Gamma) \FM(k_2,N, \epsilon_2;\Gamma) \subseteq \FM(k_1+k_2,N, \epsilon_1 \epsilon_2;\Gamma)$.  
The subsets $\FMcusppm(N)=  \oplus_{k,\epsilon} \FMcusp(k,N,\epsilon)$, and 
$\FMcusppm(N;\Gamma)=\FMpm(N;\Gamma) \cap \FMcusppm(N)$ are graded ideals 
in $\FMpm(N)$ and $\FMpm(N;\Gamma)$, respectively.  
The Fourier-Jacobi expansion map $\FJ: M_{\pm}\left( \KN \right)   \to \FMpm(N;\Gamma)$ 
is an injective homomorphism of graded rings 
that sends $S_{\pm}\left( \KN \right)$ to $\FMcusppm(N;\Gamma)$, and the inverse image of~$\FMcusppm(N;\Gamma)$ 
is $S_{\pm}\left( \KN \right)$.  
\end{lemma}
\begin{proof}
The proof is very similar to that of Lemma~\ref{lemringone}.  
\end{proof}

\begin{lemma}
\label{lemfsjfdomain}
Let $N\in \N$.  
The graded ring of formal series of Jacobi forms, 
$\FM(N)=\oplus_{k=0}^{\infty} \FM(k,N)$, is an integral domain.  
\end{lemma}
\begin{proof}
An element in $\FM(N)$ is zero precisely when each graded piece is zero.  
If a product of nonzero elements from $\FM(N)$ is zero then the product of the nonzero graded pieces of 
highest weight in each factor must be zero.  Each nonzero element in each graded piece 
has a leading term~$\phi_m \in \JkNm \subseteq {\mathcal O}(\Half_1 \times \C)$.  
The product of these leading terms cannot be zero because the ring 
${\mathcal O}(\Half_1 \times \C)$ is an integral domain.  
\end{proof}
%\smallskip 

%%%%%%%%%%%%%%%%%%%%%%%%%%%
%%%%%%%%%%%%%%%%%%%%%%%%%%%
\section{Vanishing Order of Jacobi Forms.}
\label{secvanishing}
%%%%%%%%%%%%%%%%%%%%%%%%%%%
%%%%%%%%%%%%%%%%%%%%%%%%%%%

In Corollary~{3.1} of~\cite{MR4462442}, Aoki significantly improved the known bounds on 
vanishing orders of Jacobi forms~\cite{ez85,MR1639531}, 
which were $O(m)$ in the index~$m$ for a fixed weight~$k$.    
He proved 
that a nontrivial Jacobi form of weight~$k$ and index~$m$ cannot have a vanishing order~$\nu$ 
greater than $\frac{k+1}{6}\left( \sqrt{2m+1} +1 \right)$.  As a consequence the 
$\C$-vector spaces of formal series~$\FM(k,N,\epsilon)$ are finite dimensional.  
In Theorem~\ref{thaokimain}  
we state Aoki's main theorem from~\cite{MR4462442} in the case of even weights 
and use it to improve the estimate concerning when $ J_{k,Nm}(m)=\{0\}$ enough to 
get the desired bound on the growth in~$k$ 
for dimensions of spaces of formal series of Jacobi forms,  namely, 
$\dim \FM(k,N,\epsilon) \in O(N^3k^3)$ for $N,k \in \N$.   
\smallskip

\begin{defn}
For $j\in\N$ we define $\psi_j:\N\to\Q^{+}$ by 
$$
\psi_j(u)= 
\begin{cases}  
1 
&\text{if $u=1$,}\\
{\prod}_{\substack{ p \mid u:\,\,\text{\rm $p$ prime } }}
\left( 1 - \frac{1}{p^j}\right) 
&\text{if $u \ge 2$.}  
\end{cases}
$$
\end{defn}
It is clear that $a \mid t$ implies $\psi_j(a) \ge \psi_j(t)$.  

\begin{defn}
\label{defnpsi}
We define $\psi:\N\setminus\{1\}\to\Q^{+}$ by   
$$
\psi(t)= 
\begin{cases}  
t^2\psi_2(t)=3
&\text{if $t=2$,}\\
\frac12 t^2 \psi_2(t)=\frac12 t^2 
{\prod}_{\substack{ p \mid t:\,\,\text{\rm $p$ prime } }}
\left( 1 - \frac{1}{p^2}\right) 
&\text{if $t \ge 3$.}  
\end{cases}
$$
\end{defn}
It is easy to check that $\psi(t) \ge 3$ for $ t\ge 2$.  
%\vskip0.3in
\smallskip

\begin{theorem}[\cite{MR4462442}] 
\label{thaokimain}
Let $k \in 2\N, m \in \N$.  
Let $\phi \in J_{k,m}$ have $\ord(\phi)=\mu$.  
We have 
$$
\min\left( m, m-6\mu+\frac{k}{2}\right) \ge \sum_{t: (\spadesuit)}\psi(t),
$$
where $t$ runs over all natural numbers satisfying the condition 
$$
t \ne 1 \text{ and }
\sum_{c=0}^{t-1} \psi_1\left( \gcd(t,c) \right) 
\max\left( \mu- \frac{mc(t-c)}{t^2},0\right) 
> \frac{kt}{12} \psi_2(t).  
\,\,(\spadesuit)
$$
\end{theorem}

\begin{lemma}
\label{lemaokisecond}
Let $t, \mu, m \in \N$.  
If $\mu< m$, we have 
$$
\sum_{c=0}^{t-1} \max\left( \mu- \frac{mc(t-c)}{t^2},0\right) > \frac{\mu^2 t}{2m}.  
$$  
\end{lemma}
\begin{proof} 
The condition $\mu< m$ is used in the first equality.  
\begin{align*}
&\sum_{c=0}^{t-1} \max\left( \mu- \frac{mc(t-c)}{t^2},0\right)
\ge  
\sum_{c=0}^{t-1} \max\left( \mu- \frac{mc}{t},0\right)  \\
= 
&\sum_{c=0}^{\floor(\frac{\mu t}{m})} \left( \mu- \frac{mc}{t}\right) 
=
\left( \mu
-\frac{m}{2t} \floor(\frac{\mu t}{m})\right)\left(\floor(\frac{\mu t}{m})+1  \right) \\
\ge & 
\left(\mu-\frac{\mu}{2}\right) \left(\floor(\frac{\mu t}{m})+1  \right)  
= 
\frac{\mu}{2} \left(\floor(\frac{\mu t}{m})+1  \right) > \frac{\mu^2 t}{2m}.  
\end{align*}
\end{proof}

\begin{prop}
\label{propSQuaREvanishing} 
Let $k,N,\nu \in \N$.    
If $\nu > \frac16 Nk$ then $J_{k,N\nu}(\nu)=\{0\}$.  
\end{prop}
\begin{proof}
The case of odd~$k$ follows from that of even~$k$.  
If $\phi \in J_{k,N\nu}(\nu)$ then $\phi^2 \in J_{2k,N(2\nu)}(2\nu)$. 
Since $\nu > \frac16 Nk$ we have $2\nu > \frac16 N(2k)$ and so $\phi^2=0$ 
assuming the result for even weights;  
hence $\phi=0$. 

Suppose that there is a nontrivial $\varphi \in J_{k,N\nu}(\nu)$ with 
vanishing order  
$\mu= \ord\varphi \ge \nu > \frac16 Nk$ and $k$ even.   
We will obtain a contradiction to Theorem~\ref{thaokimain}.  
For $m=N\nu$, we will contradict 
$$
%(\spadesuit)\qquad 
\min\left(m, m-6\mu+ \frac{k}{2} \right)  \ge 
\sum_{t:(\spadesuit)} \psi(t).
$$
Since $\psi$ has a positive minimum,  
%of~$3$, 
it suffices to show that the set of positive integers satisfying $(\spadesuit)$ 
is infinite.  
\smallskip

We show that all sufficiently large primes~$p$ satisfy
$$
(\spadesuit)\qquad 
\sum_{c=0}^{p-1} 
\psi_1\left( \gcd(p,c) \right)
\max\left(\mu- \frac{mc(p-c)}{p^2}, 0 \right)  >
\frac{kp}{12} \psi_2(p).
$$
We know that  
$a \mid t$ implies $\psi_1(a) \ge \psi_1(t)$.  
Since $\gcd(p,c) \mid p$, we have $\psi_1(\gcd(p,c)) \ge \psi_1(p)$.  
Therefore 
\begin{align*}
\sum_{c=0}^{p-1} 
&\psi_1\left( \gcd(p,c) \right)
\max\left(\mu- \frac{mc(p-c)}{p^2}, 0 \right)   \\
\ge 
&\psi_1(p)
\sum_{c=0}^{p-1} 
\max\left(\mu- \frac{mc(p-c)}{p^2}, 0 \right).  
\end{align*}
Lemma~\ref{lemaokisecond} says that for $t,m,\mu\in\N$ 
with $\mu<m$ we have  
$$
\sum_{c=0}^{t-1} 
\max\left(\mu- \frac{mc(t-c)}{t^2}, 0 \right) >
\frac{\mu^2t}{2m}.  
$$
We use the linear bound  $\ord \varphi \le \frac{k+2m}{12}$ of~\cite{MR1639531},  
Proposition~{3.2}, 
to 
check the hypothesis $\mu < m$:    
$$
\mu= \ord \varphi \le \frac{k+2m}{12}
< \frac{\frac{6\nu}{N}+2m}{12}
=\frac{\frac{6}{N^2}+2}{12}\, m <m.  
$$
Thus by Lemma~\ref{lemaokisecond} we have 
\begin{align*}
\sum_{c=0}^{p-1} 
&\psi_1\left( \gcd(p,c) \right)
\max\left(\mu- \frac{mc(p-c)}{p^2}, 0 \right)   \\
>
&\psi_1(p)
\frac{\mu^2p}{2m} \ge  \psi_1(p)\frac{\nu^2p}{2m} = \psi_1(p)\frac{\nu p}{2N}. 
\end{align*}
Thus a sufficient condition for a prime~$p$ to satisfy $(\spadesuit)$ 
is 
$$
\psi_1(p)\frac{\nu p}{2N} > \frac{kp}{12} \psi_2(p).  
$$
This simplifies to $\nu > \frac{Nk}{6}\left( 1+\frac{1}{p}\right)$, or 
$\nu - \frac16 Nk > \frac{Nk}{6}\frac{1}{p}$, 
which is true for all sufficiently large $p$ 
because $\nu - \frac16 Nk$ is positive.  
\end{proof}

\begin{cor}
\label{corgrowth}
For $k,N \in \N$ we have 
$\dim \FM(k,N,\epsilon) \in O(N^3k^3)$.  
\end{cor}
\begin{proof}
We have 
$\dim \FM(k,N,\epsilon) \le \sum_{j=0}^{\infty} \dim J_{k,Nj}(j+\delta)$,  
where $\delta$ is~$0$ or~$1$ as $(-1)^k \epsilon$ is~$1$ or~$-1$, 
by Lemma~{3.2} of~\cite{ipy13}.  By Proposition~\ref{propSQuaREvanishing} we 
may cap the summation at~${\floor(\frac16 Nk)}$ and, by Theorem~{2.3} 
in~\cite{ez85}, the codimension of $\Jkm^{\text{\rm cusp}}$ in $\Jkm$ 
is at most $\floor(b/2)+1$ where~$b$ is the largest integer with $b^2 \mid m$.  
Therefore   
\begin{align*}
\dim \FM(k,N,\epsilon) 
&\le \sum_{j=0}^{\floor(\frac16 Nk)} \dim J_{k,Nj}(j) 
\le \sum_{j=0}^{\floor(\frac16 Nk)} \dim J_{k,Nj}\\
&\le\dim J_{k,0} + \sum_{j=1}^{\floor(\frac16 Nk)} 
\dim J_{k,Nj}^{\text{\rm cusp}}+\tfrac12 \sqrt{Nj}+1 \\
&\le \tfrac{k+12}{12} 
+\tfrac23\left(\tfrac16 Nk+1 \right)^{3/2} +\tfrac16 Nk 
+ \sum_{j=1}^{\floor(\frac16 Nk)} \dim J_{k,Nj}^{\text{\rm cusp}}.  
\end{align*}
It is known that for $k,m\in \N$, 
$\dim J_{k,m}^{\text{\rm cusp}} \in O(km)$.  An easy estimate from the 
dimension formula~\cite{ez85} is 
$\dim J_{k,m}^{\text{\rm cusp}} \le \frac{m+1}{24}( 2k+35)$.  
From $\dim J_{k,Nj}^{\text{\rm cusp}} \in O(kNj)$ we have 
$\sum_{j=1}^{\floor(\frac16 Nk)} \dim J_{k,Nj}^{\text{\rm cusp}} \in O(N^3k^3)$.  
\end{proof}

\begin{prop}
\label{propalgebraic}
%Let $N \in \N$.  
The integral domain $\FM(N,+)=\oplus_{k=0}^{\infty} \FM(k,N,+)$ is algebraic over its subring 
$\FJ\left( M(\KNp) \right)=\oplus_{k=0}^{\infty} \FJ\left( M_k(\KNp)  \right)$.  
For $k \in \N_0$, each $\fracf \in \FM(k,N,+)$ satisfies a polynomial relation of the type 
$$
\FJ(f_0) \fracf^d + \FJ(f_1)\fracf^{d-1} + \cdots +\FJ(f_j)\fracf^{d-j} + \cdots + \FJ(f_d)=0,
$$
for some $d \in \N$, $k_0\in \N_0$, and some 
$f_j \in M_{k_0+jk}(\KN^{+})$ with $f_0$ not identically zero.       
\end{prop}
\begin{proof}
It suffices to prove the second statement.  
The group~$\Gamma=\KNp$ is commensurable with $\Sp(4,\Z)$, 
so by Theorem~{6.11} in Freitag~\cite{freitag83}, the homogeneous quotient 
field of~$M(\Gamma)$, ${\mathcal K}(\Gamma)$, has transcendence degree three.  
Take three meromorphic functions $f_1,f_2,f_3 \in {\mathcal K}(\Gamma)$ that are algebraically 
independent over~$\C$ and select a common denominator~$D \in M_{k^*}\left( \Gamma \right)$.  
We must have $k^*\in\N$ because the $f_j$ are not constant.  
We obtain four paramodular forms $g_1=Df_1$, $g_2=Df_2$, $g_3=Df_3$, and~$g_4=D$ 
in $M_{k^*}\left( \Gamma \right)$ that are algebraically independent over~$\C$.  
This follows because we may reduce to the case where any putative polynomial 
relation is homogeneous in the~$g_j$.  
\smallskip

Take $\fracf \in \FM(k,N,+)$.  
We may assume that~$\fracf$ is nontrivial because otherwise~$\fracf^1=0$ 
satisfies the conclusion with~$k_0=0$, $f_0=1\in M_0(\Gamma)$,  and $f_1=0\in M_k(\Gamma)$.   
Similarly, if $k=0$ then $\fracf=c$ is constant and~$\fracf^1-c=0$ 
satisfies the conclusion with~$k_0=0$, $f_0=1\in M_0(\Gamma)$,  and $f_1=-c\in M_0(\Gamma)$. 
So we assume $k>0$ as well.  
The four formal series~$\FJ(g_j^k) \in \FM(k^*k,N,+)$ are algebraically independent over~$\C$ 
because~$\FJ$ is a monomorphism.  Consider the list of five formal series 
$\FJ(g_j^k)$, $\fracf^{k^*} \in \FM(k^*k,N,+)$.  For any $\mu \in \N$, 
there are $\binom{\mu+4}{4}$ distinct monomials 
$x_1^{i_1}\cdots x_5^{i_5}$ in five variables with $i_1+\cdots +i_5=\mu$.  
By substitution of the five formal series into these monomials, 
we have~$\binom{\mu+4}{4}$ elements  
$\FJ(g_1^{k i_1})\cdots \FJ(g_4^{k i_4}) \fracf^{k^* i_5} \in \FM(\mu k^* k,N,+)$.  
By Corollary~\ref{corgrowth}, however, we have $\dim \FM(\mu k^* k,N,+) \in O((N k^* k)^3\mu^3)$, 
so for sufficiently large~$\mu$ there is a nontrivial $\C$-linear dependence relation 
among these~$\binom{\mu+4}{4}$ elements.   
At least one supported monomial in the dependence relation must contain a positive power of~$\fracf$ 
because the remaining four formal series are algebraically independent.  
For the same reason, when a positive power 
of~$\fracf$ is supported then its coefficient, after collecting like terms 
in powers of~$\fracf$, must be nontrivial.   
If $d \in \N$ is the highest power of~$\fracf$ that is supported in the dependence relation, 
then we may write this relation as 
$\sum_{j=0}^{d} \FJ(f_j) \fracf^{d-j}=0$ where the $f_j$ are $\C$-linear combinations of 
monomials in the four~$g_1,\ldots,g_4$.  Since $\fracf^d$ is supported, we have~$f_0$ nontrivial.  
Let the weight of~$f_0$ be~$k_0\in \N_0$ so that $f_0 \in M_{k_0}(\Gamma)$.  
The terms~$\FJ(f_j) \fracf^{d-j}$ all have the same weight, $k_0+kd$, which is the weight of~$\FJ(f_0) \fracf^{d}$.  
Therefore we may take $f_j \in M_{k_0+kj}(\Gamma)$, as required.  
\end{proof}
\vskip0.2in

%%%%%%%%%%%%%%%%%%%%%%%%%%%
%%%%%%%%%%%%%%%%%%%%%%%%%%%
\section{Invariance under subgroups of finite index in $\Gamma^0(N)$.}
\label{secinvariance}
%%%%%%%%%%%%%%%%%%%%%%%%%%%
%%%%%%%%%%%%%%%%%%%%%%%%%%% 

In Proposition~\ref{propalgebraic} of the previous section we saw that a 
formal series of Jacobi forms~$\fracf$ possessing the involution condition for~$\epsilon=+1$ 
satisfies a polynomial~$P(X)=0$ whose coefficients are formal Fourier-Jacobi expansions 
of paramodular forms.  As a consequence of this
polynomial relation we will show in this section that~$\fracf$ is invariant under a subgroup 
$\Gamma$ of finite index in $\Gamma^0(N)$.  
These arguments best take place inside the ring of formal Fourier series. 
\smallskip

The ring structure on the ring of formal Fourier series, $\ffs$, 
is defined by the Cauchy product,  
noting that for every $t \in \XtwosN$ the set $\{(t_1,t_2)\in \XtwosN \times \XtwosN: t_1+t_2=t\}$ 
is finite.  
We accordingly  use a place-holding variable~$q$ to 
write an element $\psi \in \ffs$ as 
$$
\psi= \sum_{t \in \XtwosN} a(t;\psi) q^t.  
$$
The ring of formal Fourier series, $\ffs$ is also an integral domain, which will be proven 
in Corollary~\ref{cordomainnew}.  
The Fourier expansion of a paramodular form defines a map 
\begin{align*}
\FS: M_k\left( \KN  \right) &\to \ffs \\ 
f &\mapsto \sum_{t \in \XtwosN} a(t;f) q^t.  
\end{align*}
\smallskip

Given a formal series of Jacobi forms $\fracf \in \FM(k,N)$, 
we may define the {\em associated formal Fourier series,\/} $\AFS(\fracf)\in\ffs$, by  
\begin{align*}
\AFS: \FM(k,N) &\to \ffs \\
\fracf=\sum_{m=0}^{\infty} \phi_m \xi^{Nm} &\mapsto 
\sum_{t={\tiny\nrNm} \in \XtwosN} c(n,r;\phi_m) q^t 
= \sum_{t \in \XtwosN} a(t;\fracf) q^t.     
\end{align*}
Extending by linearity we have a map $\AFS: \FM(N) \to \ffs$ 
from the ring~$\FM(N)=\oplus_{k \in \N_0} \FM(k,N)$.    
That~$\AFS$ is a ring homomorphism is an exercise using the Cauchy product 
similar to the computations demonstrating equation~{(\ref{eqcauchyprod})}.    
For $f\in  M_k\left( \KN \right)$ we have the compatibility 
$\AFS\left(\FJ(f) \right)=\FS(f)$.  
\smallskip

Formal Fourier series share some properties with formal power series~$\C[[a,b,c]]$  in three variables 
due to the following monomorphism that sends $q^t$ to the monomial 
$a^{\langle s_1,t \rangle} b^{\langle s_2,t \rangle} c^{\langle s_3,t \rangle}$.   
\begin{lemma}
\label{lemringmononew}
Set $s_1=\smallmat1000$, $s_2=\smallmat1111$, $s_3=\smallmat0001$.   
The map
\begin{align*}
\ringmono: \ffs &\to \C[[a,b,c]] \\
\psi= \sum_{t \in \XtwosN} a(t;\psi) q^t &\mapsto 
\sum_{i,j,k \in \N_0}  
\left( \sum_{\substack{t \in \XtwosN: \\ \langle s_1,t \rangle =i, \langle s_2,t \rangle =j, \langle s_3,t \rangle =k}}a(t;\psi) \right) 
a^i b^j c^k
\end{align*}
is a ring homomorphism satisfying the following properties.  
\begin{enumerate}
\item 
$\ringmono(\psi)= \sum_{i,j,k \in \N_0} 
a\left( \smallmat{i}{\frac{j-i-k}{2}}{\frac{j-i-k}{2}}{k};\psi \right)
a^i b^j c^k$,  
\item 
$\ringmono(\psi)= \sum_{t \in \XtwosN} 
a\left( t;\psi \right)
a^{\langle s_1,t \rangle} b^{\langle s_2,t \rangle} c^{\langle s_3,t \rangle}$,  
\item 
$\ringmono$ is injective.  
\end{enumerate}
\end{lemma}
\begin{proof}
Setting $t = \nrNm \in \XtwosN$ and solving the system of linear equations 
$\langle s_1,t \rangle =i$, $\langle s_2,t \rangle =j$, and $\langle s_3,t \rangle =k$, 
the unique solution is $n=i$, $r=\frac{j-i-k}{2}$, and $Nm=k$.  This proves formula~$(1)$ if we 
understand that $a(t;\psi)=0$ for $t \not\in \XtwosN$.  
Formula~$(2)$ follows from~$(1)$ since rearrangements of formal power series are equal.   
The ring homomorphism property is then formal because the $\langle s_{\ell},t \rangle$ are linear in~$t$.  
The injectivity~$(3)$ follows from formula~$(1)$.    
\end{proof}

\begin{cor}
\label{cordomainnew}
The ring $\ffs$ is an integral domain.  
\end{cor}
\begin{proof}
A ring with a monomorphism to an integral domain is an integral domain. 
The monomorphism here is $\ringmono: \ffs \to \C[[a,b,c]] $. 
\end{proof}
\smallskip

There is a copy of $\Gamma^0(N)$ inside the group of automorphisms of formal Fourier series.  
\begin{lemma}
\label{lemmajj}
For $\sigma \in \Gamma^0(N)$ define 
\begin{align*}
\jj(\sigma) : \ffs &\to \ffs \\ 
 \sum_{t\in \XtwosN} a(t;\psi) q^t & \mapsto \sum_{t\in \XtwosN} a(t[\sigma];\psi) q^t.  
\end{align*} 
The map~$\jj(\sigma)$ is an automorphism of~$\ffs$ and the map 
\begin{align*}
\jj :\Gamma^0(N) &\to \Aut\left(  \ffs  \right) \\
\sigma &\mapsto \jj(\sigma) 
\end{align*} 
is a homomorphism.  
For $\fracf \in \FM(k,N)$, we have  $\jj(\sigma) \AFS(\fracf)=\AFS(\fracf)$ if and only if 
$a(t[\sigma];\fracf)=a(t;\fracf)$ for all $t \in \XtwosN$.  
For $f \in M_k \left( \KN \right)$ we have $\jj(\sigma) \FS(f) = \FS(f)$ for all $\sigma \in \Gamma^0(N)$.  
\end{lemma}
\begin{proof}
We show $\jj(\sigma)$ is an automorphism.  
The map $\jj(\sigma)$ has $\jj(\sigma\inv)$ as an inverse and the additivity of~$\jj(\sigma)$ is clear 
so that it suffices to prove $\jj(\sigma) \left( \psi_1 \psi_2 \right) = 
\left( \jj(\sigma)\psi_1 \right) \left( \jj(\sigma)\psi_2 \right)$.  We have 
\begin{align*}
\jj(\sigma) \left( \psi_1 \psi_2 \right) &= 
\jj(\sigma)
\sum_{t\in \XtwosN} \left( \sum_{t_1,t_2\in \XtwosN:\, t_1+t_2=t}  a(t_1;\psi_1) a(t_2;\psi_2) \right) q^t 
 \\ 
&= 
\sum_{t\in \XtwosN} \left( \sum_{t_1,t_2\in \XtwosN:\, t_1+t_2=t[\sigma]}  a(t_1;\psi_1) a(t_2;\psi_2) \right) q^t .
\end{align*}
We use the equality 
$\{ (t_1,t_2) \in \XtwosN \times \XtwosN: t_1+t_2=t[\sigma]\} = 
\{ (s_1[\sigma],s_2[\sigma]) \in \XtwosN \times \XtwosN: s_1+s_2=t\}$ 
to change the index of summation.  
\begin{align*}
\jj(\sigma) \left( \psi_1 \psi_2 \right) &= 
\sum_{t\in \XtwosN} \left( \sum_{s_1,s_2\in \XtwosN:\, s_1+s_2=t} a(s_1[\sigma];\psi_1) a(s_2[\sigma];\psi_2) \right) q^t
 \\ 
&= 
\sum_{s_1\in \XtwosN}  a(s_1[\sigma];\psi_1) q^{s_1}  
\sum_{s_2\in \XtwosN}  a(s_2[\sigma];\psi_2) q^{s_2} \\
&= \left( \jj(\sigma)\psi_1 \right) \left( \jj(\sigma)\psi_2 \right).  
\end{align*}

We show that $\jj :\Gamma^0(N) \to \Aut\left(  \ffs  \right)$ is a homomorphism.  
Take $\sigma_1,\sigma_2 \in \Gamma^0(N)$ and 
$\psi =\sum_{t\in\XtwosN} a(t;\psi) q^t \in \ffs$.  
We have 
\begin{align*}
\jj(\sigma_2) \psi &= 
\sum_{t} a(t[\sigma_2];\psi) q^t
 \\  
\left( \jj(\sigma_1) \jj(\sigma_2)   \right)\psi =
\jj(\sigma_1) 
\left(  \jj(\sigma_2) \psi  \right) 
&= 
\sum_{t} a(t[\sigma_1][\sigma_2];\psi) q^t 
=\jj(\sigma_1 \sigma_2) \psi.  
\end{align*}

For $\fracf = \sum_{m=0}^{\infty} \phi_m \xi^{Nm} \in \FM(k,N)$ we have 
$\AFS(\fracf)= \sum_{t } a(t; \fracf) q^t$ 
where $a(t; \fracf) =c(n,r;\phi_m)$ for $t=\nrNm \in \XtwosN$. 
We then have 
$\jj(\sigma)\AFS(\fracf)= \sum_{t } a(t[\sigma]; \fracf) q^t$ and, 
by definition of formal series, this equals $\AFS(\fracf)$ if and only if 
$a(t[\sigma]; \fracf)=a(t; \fracf)$ for all $t \in \XtwosN$. 
The final assertion, $\jj(\sigma) \FS(f) = \FS(f)$ for all $\sigma \in \Gamma^0(N)$, 
follows from the $\Gamma^0(N)$-symmetries of $f \in M_k \left( \KN \right)$ in equation~{(\ref{eqgammaoNsym})}.  %%  
\end{proof}

The next proposition shows that a formal series of Jacobi forms satisfying 
the involution condition necessarily has additional symmetries.  

\begin{prop}
\label{propfiniteindex} 
Let $\fracf \in \FM(k,N,+)$.  
There is a subgroup $\Gamma$ of finite index in~$\Gamma^0(N)$ such that 
$\fracf \in \FM(k,N,+\,; \Gamma)$.  
\end{prop}
\begin{proof}
By Proposition~\ref{propalgebraic}, 
$\fracf$ satisfies a polynomial relation of the type 
$$
\FJ(f_0) \fracf^d + \cdots +\FJ(f_j)\fracf^{d-j} + \cdots + \FJ(f_d)=0,
$$
for some $d \in \N$, $k_0\in \N_0$, and some 
$f_j \in M_{k_0+jk}(\KN^{+})$ with $f_0$ not identically zero.    
Apply the monomorphism $\AFS: \FM(N) \to \ffs$ to obtain 
$$
\FS(f_0) (\AFS(\fracf))^d  + \cdots +\FS(f_j)(\AFS(\fracf))^{d-j} + \cdots + \FS(f_d)=0. 
$$ 
The polynomial $\sum_{j=0}^{d} \FS(f_j) X^{d-j} \in \ffs[X]$ has a root~$\AFS(\fracf)$.   
Noting by Lemma~\ref{lemmajj} that $\jj(\sigma) \FS(f_j)=\FS(f_j)$ for every $\sigma \in \Gamma^0(N)$, 
each element of the orbit $\jj\left( \Gamma^0(N) \right) \AFS(\fracf)$ is a root.  
However, the ring~$\ffs$ is an integral domain and a polynomial 
of positive degree~$d$ over an integral domain 
has at most~$d$ roots.  If we want to be definite we can find $u_1,\ldots,u_{d_1} \in \Gamma^0(N)$ 
with $d_1 \le d$ 
such that $\jj\left( \Gamma^0(N) \right) \AFS(\fracf)=\{\jj(u_1)\AFS(\fracf),\ldots,\jj(u_{d_1})\AFS(\fracf)\} $.  
The natural homomorphism~$\rho$ from $\Gamma^0(N)$ to permutations of the orbit 
$\jj\left( \Gamma^0(N) \right) \AFS(\fracf)$ is specified in this labeling by 
$\jj(\sigma) \jj(u_i) \AFS(\fracf) = \jj(u_{\rho(\sigma)i})\AFS(\fracf)$ for $i=1, \ldots, d_1$.  
Let $\Gamma=\ker(\rho)$ be the kernel of~$\rho:\Gamma^0(N)\to S_{d_1}$.    
Then $\Gamma$ is a normal subgroup of $\Gamma^0(N)$ 
of index at most~$d_1!$.  For every $\sigma \in \Gamma$ we have 
$\jj(\sigma)\AFS(\fracf)=\AFS(\fracf)$ so that, by Lemma~\ref{lemmajj}, 
$a(t[\sigma];\fracf)=a(t;\fracf)$ for all $t \in \XtwosN$.  
Thus $\fracf \in \FM(k,N,+) \cap \FM(k,N; \Gamma) = \FM(k,N,+\,; \Gamma)$ 
as claimed.  
\end{proof}

%%%%%%%%%%%%%%%%%%%%%%%%%%%
%%%%%%%%%%%%%%%%%%%%%%%%%%%
\section{Specialization.}
\label{secspecialization}
%%%%%%%%%%%%%%%%%%%%%%%%%%%
%%%%%%%%%%%%%%%%%%%%%%%%%%%

Let $\Gamma \subseteq \Gamma^0(N)$ be a subgroup of finite index.  
A formal series of Jacobi forms 
$\fracf = \sum_{m=1}^{\infty} \phi_m \xi^{Nm} \in \FMcusp(k,N;\Gamma)$ has the defining 
symmetries $a(t[\sigma];\fracf)=a(t;\fracf)$ for all $\sigma \in \Gamma$ and all $t \in \XtwosN$.  
We will construct a dense subset $\WW_1(\Gamma)$ of $\Half_1 \times \C$ where the 
formal series~$\fracf$ specializes to a holomorphic function of one variable.  More precisely, 
for each $(\tau_1,z_1) \in \WW_1(\Gamma)$, the series 
$\sum_{m=1}^{\infty} \phi_m(\tau_1,z_1) 
e\left( Nm\omega\right)$ 
will converge to a holomorphic function $\Hatz(\tau_1,z_1,\fracf)(\omega)$ on the 
neighborhood of infinity 
$\{ \omega \in \Half_1: \smallmat{\tau_1}{z_1}{z_1}{\omega} \in \Half_2\}$.  
The formal series~$\fracf$ thus converges on a dense subset of~$\Half_2$.  
\smallskip 

\begin{defn}
For $x \in \Q$, let $\Denom(x)=\min \{ n \in \N: nx \in \Z\}$ 
be the minimal positive denominator of~$x$.  
\end{defn}

\begin{defn}
\label{defSN}
Let $\Gamma \subseteq \SL(2,\Z)$ be a subgroup of finite index.    
Let $U = \smallmat0110 \in \GL(2,\Z)$.  Define 
\begin{align*}
\WW_0(\Gamma) &= 
\{ (x,y) \in \Q^2: x \in \text{$(U\Gamma U)$-orbit$(\infty)$ and } y\Denom(x)  \in \Z\},  \\
\WW_1(\Gamma) &= \{ (\tau, z) \in \Half_1 \times \C: \exists (x,y) \in \WW_0(\Gamma): z=x\tau+y\}, \\
\WW(\Gamma) &= \{ \tzzw\in \Half_2:  (\tau,z) \in \WW_1(\Gamma) \}.  
\end{align*}
\end{defn}

%%%  Fuchsian

We follow Lehner~\cite{Lehner1966} for the theory of Fuchsian groups. 
We view the Riemann sphere $\Pj^1(\C) = \C \cup \{ \infty\}$ as the extended complex plane.  
The Riemann sphere is the disjoint union 
$\Pj^1(\C) = \Pj^1(\R) \amalg \Half_1 \amalg {\overline{\Half}_1} $ 
of the extended real numbers, $\Pj^1(\R) = \R \cup \{ \infty\}$, 
and the upper and lower half planes.  
The groups $\SL(2,\R)$ and $\PSL(2,\R)=\SL(2,\R)/ \{ \pm I \}$ act on $\Pj^1(\C)$ 
by M{\"o}bius transformations and preserve this disjoint union.  
We will only consider subgroups $\Gamma \subseteq \SL(2,\R)$ and the corresponding 
transformation groups ${\bar\Gamma}= \langle \Gamma, -I \rangle/ \{ \pm I \}\subseteq \PSL(2,\R)$. 
\begin{defn}
\label{deflimitset} 
Let $\Gamma \subseteq \SL(2,\R)$.  
The {\em limit set\/} $\Lambda(\Gamma)$ of~$\Gamma$ is the set of 
$z \in \Pj^1(\C)$ such that there exists a $w \in \Pj^1(\C)$ and a sequence of 
distinct $\gamma_n \in \Gamma$ with $\lim_{n \to +\infty} \gamma_n \langle w \rangle =z$.  
\end{defn}
For subgroups $\Gamma_1, \Gamma_2 \subseteq \SL(2,\R)$, 
if $\Gamma_1$ has finite index in $\Gamma_2$ then $\Lambda(\Gamma_1)=\Lambda(\Gamma_2)$; 
this is the theorem in section~{2C} of~\cite{Lehner1966}, page~{11}.    
The corollary in section~{2F}, page~{14}, is that either $\Lambda(\Gamma)=\Pj^1(\C)$ or 
$\Lambda(\Gamma) \subseteq\Pj^1(\R)$.  We mention some terminology to assist readers who 
use a different reference.  
A group $\bar\Gamma \subseteq \PSL(2,\R)$ is {\em Fuchsian\/} when 
$\Lambda(\Gamma) \subseteq\Pj^1(\R)$.  
The theorem in section~{2F} of~\cite{Lehner1966}, page~{13}, 
characterizes Fuchsian groups as the discrete subgroups of~$\PSL(2,\R)$.  
For us, the salient result is Theorem~{3} in section~{3E}, page~{21}.  
\begin{theorem}[{\cite{Lehner1966}}]
\label{theoremLehner}
Let $\Gamma \subseteq \SL(2,\R)$ be a subgroup.  
If $ F \subseteq \Pj^1(\C)$ is a closed set containing at least two points, such that 
$\Gamma F \subseteq F$, then $F \supseteq \Lambda(\Gamma)$.  
\end{theorem}
As Lehner comments on page~{21}, this theorem may be rephrased: 
when $\Lambda(\Gamma)$ has more than one point, $\Lambda(\Gamma)$ 
is the smallest closed $\Gamma$-invariant set containing at least two points.  
\smallskip

An example of a Fuchsian group is~$\PSL(2,\Z)$.  
The orbit of~$\infty$ is $\Pj^1(\Q)$ and $\Lambda(\SL(2,\Z))=\Pj^1(\R)$.  
This implies that any subgroup $\Gamma \subseteq \SL(2,\Z)$ of finite index 
also has $\Lambda(\Gamma)=\Pj^1(\R)$.  
For a ring $R \subseteq \R$, set 
$\Ponezero(R)=\{ \smallmat{a}bcd \in \SL(2,R): c=0 \}$.  
The stabilizer of~$\infty$ in~$\SL(2,\R)$ is~$\Ponezero(\R)$  
and we have an orbit-stabilizer bijection 
$\SL(2,\Z)/ \Ponezero(\Z) \leftrightarrow \Pj^1(\Q)$ given by 
sending $\gamma \Ponezero(\Z) \mapsto \gamma \langle \infty \rangle$.  
This bijection shows that a subgroup of finite index in~$\SL(2,\Z)$ cannot 
stabilize~$\infty$.  
We require the following ergodic corollary of Theorem~{\ref{theoremLehner}}.  

\begin{lemma} 
\label{lemmamosher}
Let $\Gamma \subseteq \SL(2,\Z)$ be a subgroup of finite index.   
The $\Gamma$-orbit of~$\infty$ is dense in $\Lambda(\Gamma)= \Pj^1(\R)$.  
\end{lemma}  
\begin{proof}
We have  $\Lambda(\Gamma)=\Lambda(\SL(2,\Z))= \Pj^1(\R)$ because 
$\Gamma$ has finite index in~$\SL(2,\Z)$.  Furthermore, 
$\Gamma$ cannot stabilize~$\infty$ for the same reason.  
Therefore, the $\Gamma$-orbit of~$\infty$ has at least two points, as does its closure 
$F=\overline{\text{$\Gamma$-orbit$(\infty)$}} \subseteq  \Pj^1(\R)$.  
We check that 
the closure~$F$ remains $\Gamma$-invariant.  
Take $z \in F$ and $\gamma \in \Gamma$.    
We have $z=\lim_{n \to +\infty} \gamma_n \langle \infty \rangle$ 
for some $\gamma_n \in \Gamma$, so that  
$\gamma \langle z \rangle=\lim_{n \to +\infty} (\gamma\gamma_n) \langle \infty \rangle$ 
for $\gamma\gamma_n \in \Gamma$.  
Hence $\gamma\langle z \rangle \in F$ and $F$ is $\Gamma$-invariant. 
Thus, by Theorem~\ref{theoremLehner} we have 
$F = \Lambda(\Gamma)=\Pj^1(\R)$.  The 
$\Gamma$-orbit of~$\infty$ is thus dense in $\Pj^1(\R)$. 
\end{proof} 
\smallskip

\begin{lemma}
\label{lemdensity}
Let $\Gamma \subseteq \SL(2,\Z)$ be a subgroup of finite index.  
The sets $\WW_0(\Gamma)$, $\WW_1(\Gamma)$, and $\WW(\Gamma)$ are dense in 
$\R^2$, $\Half_1 \times \C$, and~$\Half_2$, respectively.  
\end{lemma}
\begin{proof} 
It suffices to prove $\WW_0(\Gamma)$ is dense in~$\R^2$.  
Take $(\xi,\eta) \in \R^2$ and any neighborhoods $X$ of $\xi$ and $Y$ of $\eta$.  
We just need to find an element of $\WW_0(\Gamma)$ in $X \times Y$.  
Pick an irrational number $\xi_o \in X$.  

The group $U\Gamma U$ has finite index in~$\SL(2,\Z)$ and so
$\Lambda(U\Gamma U)=\Pj^1(\R)$.  
By Lemma~\ref{lemmamosher}, 
the $(U\Gamma U)$-orbit$(\infty)$ is dense in $\Lambda(U\Gamma U)=\Pj^1(\R)$, 
and hence $\Q \cap (U\Gamma U)$-orbit$(\infty)=(U\Gamma U)$-orbit$(\infty)\setminus \{\infty\}$ is dense in $\R$.  
Accordingly there is a sequence $x_j \in \Q \cap (U\Gamma U)$-orbit$(\infty)$ 
with $\lim_j x_j = {\xi_o}$, and we have $x_j \in X$ for all sufficiently large~$j$.  
Set $D_j = \Denom(x_j)$.  Since ${\xi_o}$ is irrational we have $\lim_j D_j = +\infty$.  
Set $y_j =\frac{1}{D_j} \floor(\eta D_j)$ so that $(x_j,y_j) \in \WW_0(\Gamma)$ and 
$\lim_j y_j = \eta$.  For all sufficiently large~$j$ we have $(x_j,y_j) \in X \times Y$, 
which shows that $\WW_0(\Gamma) \cap (X \times Y)$ is nonempty, and hence that 
$\WW_0(\Gamma)$ is dense in~$\R^2$.  
\end{proof}

Let $\Gamma$ be a subgroup of~$\SL(2,\Z)$ with finite index.   
For $f \in S_k\left( \Gamma \right)$, 
define the {\em Hecke bound\/} 
$\HB(f)= \sup_{\tau \in \Half_1} (\IM \tau)^{k/2} \vert f(\tau) \vert$.  
The Hecke bound is a norm on the vector space $ S_k\left( \Gamma \right)$ and the 
Fourier coefficients $a(t, \cdot): S_k\left( \Gamma \right) \to \C$, given by $f \mapsto a(t;f)$, 
 become bounded linear functionals.  
Hence the Fourier coefficients are continuous in~$f$ in the topology on 
$S_k\left( \Gamma \right)$ induced by this norm.  For Jacobi cusp forms 
$\phi \in J_{k,m}^{\text{\rm cusp}}$ the Hecke bound~\cite{ez85}, compare page~$27$, is given by 
$$
\HB(\phi)= \sup_{(\tau,z) \in \Half_1 \times \C} 
(\IM \tau )^{k/2}\, e^{-2\pi m \frac{(\IM z)^2}{\IM \tau}}\,
\vert \phi(\tau,z) \vert.  
$$

\begin{lemma}
\label{lemHB}
Let $\phi \in J_{k,m}^{\text{\rm cusp}}$ with $m \in \N$.  
For all $(n,r) \in \N \times \Z$ we have the bound 
$\vert c(n,r;\phi) \vert \le \HB(\phi)  
\left( \frac{e\pi}{mk}\vert 4mn-r^2\vert \right)^{k/2}$.  
\end{lemma}
\begin{proof} 
Let 
$\phi(\tau,z) = \sum_{n \in \N,\, r \in \Z} c(n,r; \phi)\, e(n\tau + r z)$ 
be the Fourier series with  
$\tau = u+iv \in \Half_1$ and $z=x+iy \in \C$.  
By the familiar formula for the Fourier coefficients, we have 
\begin{align*}
c(n,r;\phi)  &=  
\iint_{[0,1]^2} 
\phi(\tau,z) e(-n\tau - r z)\,du\,dx  \\  
\vert  c(n,r;\phi) \vert &\le  
\iint_{[0,1]^2} 
\vert \phi(\tau,z)\vert\,  e^{2\pi nv} e^{2\pi ry}\,du\,dx  \\  
 &\le  
\iint_{[0,1]^2} 
\HB(\phi) 
v^{-k/2}  e^{2\pi m \frac{y^2}{v}} 
e^{2\pi(nv+ry)} \,du\,dx  \\  
 &=   
\HB(\phi) 
v^{-k/2}  e^{2\pi (m \frac{y^2}{v}+ry+nv  )}.  
\end{align*}
Choosing $y=-\frac{rv}{2m}$ we have   
$m \frac{y^2}{v}+ry+nv = \frac{4mn-r^2}{4m}  v$.  
Thus we have 
$ \vert  c(n,r;\phi) \vert \le    \HB(\phi) v^{-k/2} e^{2\pi \frac{4mn-r^2}{4m}  v}$.  
For $4mn-r^2 > 0$ this is minimized by 
$v= \frac{km}{\pi(4mn-r^2)}$ proving 
$\vert c(n,r;\phi) \vert \le \HB(\phi)  
\left( \frac{e\pi}{mk}(4mn-r^2) \right)^{k/2}$.  
This proves the result when $4mn-r^2 > 0$, and   
when $4mn-r^2 \le 0$ we have $c(n,r;\phi)=0$ so the Lemma's conclusion holds as well.  
\end{proof}

The following proposition gives a specialization of a 
cuspidal formal series~$\fracf \in \FMcusp(k,N;\Gamma)$ for each $(x,y)\in \WW_0(\Gamma) $.  
Each specialization depends only upon a finite number of the Fourier-Jacobi coefficients of~$\fracf$, 
although this finite number may increase with the denominator of~$x$.    
\vskip0.1in

\begin{prop}
\label{propaokinew} 
Let $\Gamma \subseteq \Gamma^0(N)$ be a subgroup of finite index.  
Let $\fracf = \sum_{m=1}^{\infty} \phi_m \xi^{Nm} \in \FMcusp(k,N;\Gamma)$.  
For $(x,y)\in \WW_0(\Gamma) $, set $D=\Denom(x)$.  We have 
$
f_m(\tau):=e\left( Nmx^2\tau\right) 
\phi_m(\tau, x\tau+y) \in S_k\left( \Gamma(D^2) \right).  
$  
The sequence of Hecke bounds satisfies 
$$\HB(f_m) \in O\left(m^{\frac{k+1}{2}}\right),$$ 
where the implied constant depends only on~$\fracf$ and~$x$.  %%
\end{prop}
\begin{proof}
Take $\fracf = \sum_{m=1}^{\infty} \phi_m \xi^{Nm} \in \FMcusp(k,N;\Gamma)$.   
For each~$m \in \N$ there is a positive constant~$K_m$ such that 
$\vert c(n,\ell;\phi_m) \vert \le K_m \left( \frac{4Nmn-\ell^2}{4Nm} \right)^{k/2}$ 
for all $(n,\ell) \in \supp(\phi_m)$.  Here we use the assumption that each~$\phi_m$ is a 
cusp form and Lemma~\ref{lemHB}.  
By Proposition~{2.2} of~\cite{MR3799156} or Theorem~{4.2} of~\cite{MR4462442}, 
for example, we have 
$f_m \in S_k\left( \Gamma(D^2) \right)$ because $Dx,Dy \in \Z$.  
This elliptic modular cusp form has the 
Fourier expansion
$$
f_m(\tau) = \sum_{v > 0:\, D^2v \in \Z} a(v;f_m) q^v; \   
a(v;f_m)=\sum_{(n,\ell)\in \Lambda_{Nm,x,v}} e(y\ell) c(n,\ell;\phi_m), 
$$
for $\Lambda_{Nm,x,v}=
\{(n,\ell) \in \N \times \Z: Nm x^2+x\ell +n=v, \, 4Nmn> \ell^2\}$.  
Next we check $\vert \Lambda_{Nm,x,v} \vert \le 4 \sqrt{Nmv}+1$.  
For any $(n,\ell) \in  \Lambda_{Nm,x,v}$ we  
first note that~$n$ is determined by~$\ell$, and that~$\ell$ satisfies 
$$
\text{ $ (2Nmx+\ell)^2 + 4Nmn-\ell^2=4Nm(Nmx^2+\ell x +n) = 4Nmv$.}  
$$
Thus $\vert 2Nmx+\ell \vert \le 2\sqrt{Nmv}$ and the number of integers~$\ell$ in this  
interval is at most the length plus one, showing  $\vert \Lambda_{Nm,x,v} \vert \le 4 \sqrt{Nmv}+1$.  
Pick an $L \in \N$ such that the Fourier coefficients~$a(v;f)$ for $v \le L$ determine 
$f \in S_k\left( \Gamma(D^2) \right)$. 
By virtue of the valence inequality, 
the choice $L=L(k,D)=\floor\left( \frac{k}{12} [\SL(2,\Z): \Gamma(D^2)] \right)$ will work.   
The following supremum exists because the denominator is nonzero, the quotient 
is continuous in~$f$, and the supremum may 
be taken over a compact sphere.  
%The denominator in the following supremum is nonzero. 
$$
C(k,D)= \sup_{f \in S_k\left( \Gamma(D^2) \right)\setminus \{0\} }
\dfrac{\HB(f)}{\sum_{0<v \le L(k,D)} \vert a(v;f) \vert }
$$
Therefore we have 
\begin{align*}
\HB(f_m) &\le C(k,D) \sum_{0<v \le L(k,D)} \vert a(v;f_m) \vert \\
&\le C(k,D) \sum_{0<v \le L(k,D)} \sum_{(n,\ell)\in \Lambda_{Nm,x,v}} \vert c(n,\ell;\phi_m) \vert \\
&\le C(k,D)  \sum_{(n,\ell)\in {\tilde\Lambda}_{Nm,x,L}} \vert c(n,\ell;\phi_m) \vert 
\end{align*}
for ${\tilde \Lambda}_{Nm,x,L}=
\{(n,\ell) \in \N \times \Z: Nm x^2+x\ell +n \le L, \, 4Nmn> \ell^2\}$ or, 
equivalently, $ {\tilde \Lambda}_{Nm,x,L}=\cup_{0<v\le L}\,\Lambda_{Nm,x,v}$.  
Each $\Lambda_{Nm,x,v}$ has at most $4 \sqrt{Nmv}+1$ elements and there are at most 
$D^2L$ of them so we have 
$\vert {\tilde \Lambda}_{Nm,x,L}\vert \le D^2L (4 \sqrt{NmL}+1)$.  
Note that $ (n,\ell) \in {\tilde \Lambda}_{Nm,x,L}$ implies that $4Nmn-\ell^2 \le 4NmL$.  
This is because 
$$
\text{\footnotesize $4Nmn-\ell^2 \le (2Nmx+\ell)^2 + 4Nmn-\ell^2=4Nm(Nmx^2+\ell x +n) \le 4NmL$.}  
$$

We now come to the main point where we use 
the $\Gamma$-conditions for~$\fracf$.   
Our assumption $(x,y)\in \WW_0(\Gamma)$ gives us $x \in (U\Gamma U)$-orbit$(\infty)$.  
Accordingly there exists a matrix 
$\smallmat{a}{\eta}{c}{\xi} \in U\Gamma U$ with $x=a/c$.  
Letting $D=\Denom(x)$, the conditions $\gcd(a,c)=1$ and $x=a/c$ imply $c = \pm D$ and $a=cx$.  
Moving this matrix to~$\Gamma$ we define 
$$
\sigma= \smallmat{\xi}{c}{\eta}{a}= 
U \smallmat{a}{\eta}{c}{\xi} U \in \Gamma.  
$$
We note that $\Gamma \subseteq \Gamma^0(N)$ implies $N \mid D$.  
For  $\smallmat{n}{\ell/2}{\ell/2}{Nm} \in \XtwoN$, compute 
$\sigma' \smallmat{n}{\ell/2}{\ell/2}{Nm} \sigma = \smallmat{n_1}{\ell_1/2}{\ell_1/2}{Nm_1}\in \XtwoN$, where 
\begin{align*}
n_1 &= Nm \eta^2 + \ell \xi \eta + n \xi^2, \\
\ell_1&=c\left(2n\xi + \ell(\eta+\xi x) + 2Nmx\eta  \right), \\
m_1 &= \tfrac{D}{N} D\left( Nmx^2+\ell x +n \right) .  
\end{align*}
For $(n,\ell) \in {\tilde \Lambda}_{Nm,x,L}$ we have $Nmx^2+\ell x +n\le L$ so that 
$m_1 \le \frac{D^2}{N}L$, 
as well as  $\smallmat{n}{\ell/2}{\ell/2}{Nm}$ and $\smallmat{n_1}{\ell_1/2}{\ell_1/2}{Nm_1} \in \XtwoN$.   
By the $\Gamma$-conditions of equation~{(\ref{eqgammaoNsymformal})} for~$\fracf$ we see 
\begin{align*}
&\vert c(n,\ell; \phi_m) \vert = \vert c(n_1,\ell_1; \phi_{m_1}) \vert 
\le K_{m_1} 
\left( \frac{4Nm_1n_1-\ell_1^2}{4Nm_1}  \right)^{k/2} \\
&=K_{m_1}  \left( \frac{4Nmn-\ell^2}{4Nm_1}  \right)^{k/2} 
\le K_{m_1}  \left( \frac{4NmL}{4Nm_1}  \right)^{k/2} 
\le K_{m_1}  \left( Lm  \right)^{k/2}, 
\end{align*}
if we remember that $(n,\ell) \in {\tilde \Lambda}_{Nm,x,L}$ implies 
$0< 4Nmn-\ell^2 \le 4NmL$ and note $m_1 \ge 1$.  
Therefore, for $(n,\ell) \in {\tilde \Lambda}_{Nm,x,L}$, we have 
\begin{align*}
\vert c(n,\ell; \phi_m) \vert &\le K_{m_1}  \left( Lm  \right)^{k/2} 
\le \left( \max_{1 \le j \le \frac{D^2}{N}L} K_j \right) L^{k/2} m^{k/2} \\ 
&=K_{D,N,L} L^{k/2} m^{k/2}   
\end{align*}
since $m_1 \le \frac{D^2}{N}L$.  
The constant~$K_{D,N,L}$ grows with the denominator of $(x,y) \in \WW_0(\Gamma)$ 
but in each individual case controls the growth of all~$\HB(f_m)$ in terms of the 
finite number of~$\HB(\phi_j)$ with $j \le \frac{D^2}{N}L$.  
Therefore we have 
\begin{align*}
\HB(f_m)
&\le C(k,D)  \sum_{(n,\ell)\in {\tilde\Lambda}_{Nm,x,L}} \vert c(n,\ell;\phi_m) \vert \\  
&\le C(k,D)  \sum_{(n,\ell)\in {\tilde\Lambda}_{Nm,x,L}} K_{D,N,L}L^{k/2} m^{k/2} \\
&= C(k,D)\,   \vert {\tilde\Lambda}_{Nm,x,L}\vert\, K_{D,N,L}L^{k/2} m^{k/2} \\
&\le C(k,D) K_{D,N,L}L^{k/2}  D^2L (4 \sqrt{NmL}+1)  m^{k/2} \\
&\le C_2(k,D,N,L)  (4 \sqrt{NmL}+1)  m^{k/2}.  
\end{align*}
Thus $\HB(f_m) \in O\left(m^{\frac{k+1}{2}}\right)$ as was to be shown.  
\end{proof}

 Consider a formal series 
$\fracf = \sum_{m=0}^{\infty} \phi_m \xi^{Nm} \in \FM(k,N)$. 
For a point $\Omega_o=\smallmat{\tau_o}{z_o}{z_o}{\omega_o} \in \Half_2$, 
if the series $\sum_{m=0}^{\infty} \phi_m(\tau_o,z_o) e(Nm\omega_o)$ 
converges absolutely,  we say that the formal series~$\fracf$ converges absolutely at~$\Omega_o$.  

\begin{cor}
\label{coraokinew}
Let $\Gamma \subseteq \Gamma^0(N)$ be a subgroup of finite index.  
Let $\fracf = \sum_{m=1}^{\infty} \phi_m \xi^{Nm} \in \FMcusp(k,N;\Gamma)$.  
The formal series~$\fracf$ converges absolutely on $\WW(\Gamma)$.    
For $(\tau_1,z_1) \in \WW_1(\Gamma)$ set $\eta_1=(\IM {z_1})^2/{\IM {\tau_1}}$.  
The function $\Hatz(\tau_1,z_1,\fracf): \Horo(\eta_1) \to \C$ defined by 
$$
\Hatz(\tau_1,z_1,\fracf)(\omega) = \sum_{m=1}^{\infty} \phi_m(\tau_1,z_1) 
e\left( Nm\omega\right)
$$ 
is holomorphic on $\Horo(\eta_1)$.   
\end{cor}
\begin{proof}
Pick $\Omega=\smallmat{\tau}{x\tau+y}{x\tau+y}{\omega} \in \WW(\Gamma)$ 
with $(x,y) \in \WW_0(\Gamma)$.  By Proposition~\ref{propaokinew}, 
$\HB(f_m) \le C(\fracf,x) m^{\frac{k+1}{2}}$ 
for some constant~$C(\fracf,x)$ and the elliptic forms 
$f_m \in S_k\left( \Gamma(D^2) \right)$ given by 
$
f_m(\tau)=e\left( Nmx^2\tau\right) 
\phi_m(\tau, x\tau+y)  
$.  
From $(\IM \tau)^{k/2} \vert f_m(\tau) \vert \le \HB(f_m) \le C(\fracf,x) m^{\frac{k+1}{2}}$ 
we infer 
\begin{align*}
\vert  \phi_m(\tau, x\tau+y)  \vert &\le C(\fracf,x) 
(\IM \tau)^{-k/2}\, e^{2\pi Nm x^2 (\IM \tau)}
m^{\frac{k+1}{2}}, \\
\limsup_{m \to +\infty} \vert  \phi_m(\tau, x\tau+y)  \vert^{\frac1m}  &\le 
e^{2\pi N x^2 (\IM \tau)}.  
\end{align*}
The radius of convergence~$R$ of the series 
$\sum_{m \in \N} \phi_m(\tau,x\tau+y) \xi^{Nm}$ thus satisfies 
$R \ge e^{-2\pi N x^2 (\IM \tau)}$, and so this series converges absolutely for 
$\vert \xi \vert^N <  e^{-2\pi N x^2 (\IM \tau)^2}$.  
In particular the series converges at~$\Omega$ if 
$\IM\omega > x^2 (\IM \tau) = \frac{(\IM {x\tau+y})^2}{\IM \tau}$. 
The condition $\IM \omega \IM \tau > (\IM {x\tau+y})^2$, however, is just 
$\Omega=\smallmat{\tau}{x\tau+y}{x\tau+y}{\omega} \in \Half_2$.  
\smallskip

Consider $(\tau_1,z_1) \in \WW_1(\Gamma)$ so that $(\tau_1,z_1)= (\tau_1, x\tau_1+y)$ 
for some $(x,y)\in \WW_0(\Gamma)$.  
We have seen that the power series in $e(\omega)$,  
\begin{equation}
\label{eqHatzholo}
\sum_{m \in \N} \phi_m(\tau_1,z_1) e\left(Nm\omega  \right), 
\end{equation}
converges when $\smallmat{\tau_1}{z_1}{z_1}{\omega} \in \Half_2$.  
A power series that converges on an open set is holomorphic there.  
Thus  $\Hatz(\tau_1,z_1,\fracf)$ 
is holomorphic on the open set 
$\{ \omega \in \Half_1: \smallmat{\tau_1}{z_1}{z_1}{\omega} \in \Half_2\}=\Horo(\eta_1)$.  
\end{proof}
\vskip0.2in

%%%%%%%%%%%%%%%%%%%%%%%%%%%
%%%%%%%%%%%%%%%%%%%%%%%%%%%
\section{Locally Bounded.}
\label{seclocallybounded}
%%%%%%%%%%%%%%%%%%%%%%%%%%%
%%%%%%%%%%%%%%%%%%%%%%%%%%%

Locally bounded families of holomorphic functions possess remarkable 
convergence properties.  In Theorem~\ref{lemlocallybounded} we show 
that, for subgroups $\Gamma \subseteq \Gamma^0(N)$ of finite index,  
the partial sums of the formal series 
$\fracg = \sum_{m=1}^{\infty} \phi_m \xi^{Nm}\in \FMcusp(k,N,\epsilon;\Gamma)$ 
are locally bounded on~$\Half_2$ if~$\fracg$ is integral over 
$\FJ\left(\FMpm(N;\Gamma)\right)$.  
The following homomorphism respects the ring structure but forgets the grading on~$\FMpm(N)$.  
\begin{lemma}
\label{lemAtz}
For each $(\tau_1,z_1) \in \Half_1 \times \C$ there is a ring homomorphism   
$\Atz(\tau_1,z_1): \FMpm(N) \to \C[[\xi]]$ defined by sending 
$\fracf = \sum_{m=0}^{\infty} \phi_m \xi^{Nm} \in \FM(k,N,\epsilon)$ 
to $\Atz(\tau_1,z_1)\fracf = \sum_{m=0}^{\infty} \phi_m(\tau_1,z_1) \xi^{Nm} \in \C[[\xi]]$  
and extending additively to $\FMpm(N)=\oplus_{k,\epsilon} \FM(k,N,\epsilon)$.  
\end{lemma}
\begin{proof}
The ring structures are defined by the Cauchy product rule, so substitution is a homomorphism.  
\end{proof}

Let $D=\{z\in\C: \vert z \vert <r\}$ be the open disk about the origin in~$\C$ with positive radius~$r$.  
Let $\Sps(D)=\{w \in \C[[\xi]]: \text{$w$ converges on~$D$} \}$ 
be the ring of power series that converge on~$D$.  
The map $\Taylor: {\mathcal O}(D) \to \Sps(D)$ 
that sends a holomorphic function on~$D$ to its Taylor series about the origin is a ring isomomorphism.  
For $\xi_o \in D$ and $w=\sum_{n=0}^{\infty} w_n \xi^n \in \Sps(D)$, 
let $w[\xi_o]=\sum_{n=0}^{\infty} w_n \xi_o^n \in\C$ be the sum of~$w$ at~$\xi:=\xi_o$.  
The evaluation map $[\xi_o]: \Sps(D) \to \C$ sending $w$ to $w[\xi_o]$ is a ring homomorphism 
with the property that $\left( \Taylor\inv w \right)(\xi_o)=w[\xi_o]$.  

\begin{lemma}
\label{lemAtzs0}
For each $(\tau_1,z_1) \in \Half_1 \times \C$, define 
$\eta_1=\frac{(\IM {z_1})^2}{\IM {\tau_1}}$,  a radius 
$r_1=e^{-2\pi \eta_1}$, and a disk $D_1=\{\xi\in\C: \vert \xi \vert <r_1\}$.  
The ring homomorphism   
$\Atz(\tau_1,z_1): \FMpm(N) \to \C[[\xi]]$ sends 
$\FJ\left(M_{\pm}(N) \right)$ into $\Sps(D_1)$.    
For $\omega \in \Horo(\eta_1)$, 
set $\Omega_1(\omega)=\smallmat{\tau_1}{z_1}{z_1}{\omega}\in \Half_2$.  
For $f \in M_k\left( \KN \right)^{\epsilon}$ and $\omega \in \Horo(\eta_1)$, we have 
$\left(\Atz(\tau_1,z_1) \FJ(f)\right)[e(\omega)] {=}f\left(\Omega_1(\omega) \right)$.  
\smallskip

Let $\Gamma \subseteq \Gamma^0(N)$ be a subgroup of finite index.  
For $(\tau_1,z_1) \in \WW_1(\Gamma)$, 
$\Atz(\tau_1,z_1)$  sends 
%  $\Atz(\tau_1,z_1): \FMpm(N) \to \C[[\xi]]$ sends 
$\FMcusppm(N;\Gamma)$ into $\Sps(D_1)$.  
For  $\fracf \in \FMcusp(k,N,\epsilon;\Gamma)$ and $\omega \in \Horo(\eta_1)$ 
we have 
$\left( \Atz(\tau_1,z_1) \fracf \right)[e(\omega)]=\Hatz(\tau_1,z_1,\fracf)(\omega)$ 
for the holomorphic function~$\Hatz(\tau_1,z_1,\fracf)$ defined in Corollary~\ref{coraokinew}.  
%\smallskip
%\vskip0.1in
%  
The map $\Hatz(\tau_1,z_1): \FMcusppm(N;\Gamma) \to {\mathcal O}(\Horo(\eta_1))$ 
that sends~$\fracf$ to $\Hatz(\tau_1,z_1,\fracf)$ is linear and multiplicative. 
\end{lemma}
\begin{proof}  
We know that $\Atz(\tau_1,z_1)$ is a ring homomorphism from Lemma~\ref{lemAtz}.  
We need to check that the relevant power series converge in~$D_1$ and correctly 
label their values.  
For $f \in M_k\left( \KN \right)^{\epsilon}$ the Fourier-Jacobi expansion 
$f(\Omega)=\sum_{m=0}^{\infty} \phi_m(\tau,z) e(Nm\omega)$ converges for 
$\Omega=\tzzw\in \Half_2$.  The formal series is 
$\FJ(f)=\sum_{m=0}^{\infty} \phi_m \xi^{Nm}$ and the value of~$\Atz(\tau_1,z_1)$ on it is 
$\Atz(\tau_1,z_1)\FJ(f)=\sum_{m=0}^{\infty} \phi_m(\tau_1,z_1) \xi^{Nm} \in \C[[\xi]]$.  
For any $\xi_o \in D_1$, there is an $\omega_o \in \Horo(\eta_1)$ with $\xi_o=e(\omega_o)$ and 
$\Omega_1(\omega_o) \in \Half_2$. Accordingly, 
$f(\Omega_1(\omega_o))=\sum_{m=0}^{\infty} \phi_m(\tau_1,z_1) \xi_o^{Nm}$ is convergent 
and, since $\xi_o \in D_1$ was arbitrary, $\Atz(\tau_1,z_1)\FJ(f) \in \Sps(D_1)$.  
For any $\omega \in \Horo(\eta_1)$ we have $e(\omega) \in D_1$ and 
$\left(\Atz(\tau_1,z_1) \FJ(f)\right)[e(\omega)] = 
\sum_{m=0}^{\infty} \phi_m(\tau_1,z_1) e(Nm\omega)=
f\left(\Omega_1(\omega) \right)$.  
\smallskip

Assume that $(\tau_1,z_1) \in \WW_1(\Gamma)$.  
Since $\fracf \in \FMcusp(k,N,\epsilon;\Gamma) \subseteq \FMcusp(k,N;\Gamma)$, 
we can apply Corollary~\ref{coraokinew} to assure the convergence of the series 
$
\Hatz(\tau_1,z_1,\fracf)(\omega){=}
\sum_{m=1}^{\infty} \phi_m(\tau_1,z_1) e(Nm\omega) 
$ 
for $\omega \in \Horo(\eta_1)$.    
Any $\xi \in D_1$ may be written $\xi=e(\omega)$ for $\omega \in \Horo(\eta_1)$, so the power series 
$
\Atz(\tau_1,z_1)\fracf=\sum_{m=1}^{\infty} \phi_m(\tau_1,z_1) \xi^{Nm}
$ 
converges on~$D_1$ and so $\Atz(\tau_1,z_1)\fracf \in \Sps(D_1)$.  
The value of the convergent series at $e(\omega) \in D_1$ is given by 
$\left( \Atz(\tau_1,z_1) \fracf \right)[e(\omega)]=
\sum_{m=1}^{\infty} \phi_m(\tau_1,z_1) e(Nm\omega) = 
\Hatz(\tau_1,z_1,\fracf)(\omega)$.
\smallskip

 Both $\Atz(\tau_1,z_1)$ and the 
evaluation map $[e(\omega)]$ are ring homomorphisms so their composition, 
restricted to the ideal~$\FMcusppm(N;\Gamma)$ is linear and multiplicative.  
This shows, for all $\omega \in \Horo(\eta_1)$, that we have 
\begin{align*}
&\Hatz(\tau_1,z_1,\fracf_1\fracf_2)(\omega)=
\left( \Atz(\tau_1,z_1) \fracf_1\fracf_2 \right)[e(\omega)] = \\  
&\left( \Atz(\tau_1,z_1) \fracf_1 \right)[e(\omega)]\left( \Atz(\tau_1,z_1) \fracf_2 \right)[e(\omega)]
=\Hatz(\tau_1,z_1,\fracf_1)(\omega)\Hatz(\tau_1,z_1,\fracf_2)(\omega),
\end{align*}   
which is what  
$\Hatz(\tau_1,z_1,\fracf_1\fracf_2)
=\Hatz(\tau_1,z_1,\fracf_1)\Hatz(\tau_1,z_1,\fracf_2)$ means.  
\end{proof}

\begin{lemma}
\label{lempoly}
Let $d\in \N$ and $z \in \C$ be given.  
Let a monic polynomial $P\in \C[X]$ of degree~$d$  be given by
$
P(X)=
X^d+ \sum_{j=1}^d a_j X^{d-j}
$ 
for $a_1, \ldots, a_d \in \C$.  
If $P(z)=0$ then $\vert z \vert \le 1+ \sum_{j=1}^d \vert a_j \vert$.  
\end{lemma}
\begin{proof}
Apply the Triangle Inequality to the case $\vert z \vert \ge 1$.  
\end{proof}

\begin{theorem}
\label{lemlocallybounded}
Let $N,k,d \in \N$.  Let $\epsilon \in \{ \pm 1\}$.  
Let $\Gamma \subseteq \Gamma^0(N)$ be a subgroup of finite index.  
Let $\fracg = \sum_{m=1}^{\infty} \phi_m \xi^{Nm} 
\in \FMcusp(k,N,\epsilon;\Gamma)$ satisfy the monic polynomial relation 
$P(\fracg)=0$ in $\FM(kd,N,\epsilon^d;\Gamma)$ for 
$$
P(X)=X^d+
\cdots + \FJ(g_j) X^{d-j}+ \cdots + \FJ(g_d)
$$
with $g_j \in M_{kj}\left( \KN \right)^{\epsilon^j}$ for $j=1,\ldots,d$.  

Then the sequence of partial sums 
$\sum_{m=1}^{M} \phi_m(\tau,z) e\left( Nm\omega\right)$ for $M \in \N$ 
is locally bounded on $\tzzw\in \Half_2$.  
\end{theorem}
\begin{proof}
Pick $\Omega_o \in  \Half_2$. 
We need to make a neighborhood~$B_1$ of~$\Omega_o$ and a positive constant~$A_1$ 
such that for all $M \in \N$ and all $\tzzw\in B_1$ we have 
$\vert \sum_{m=1}^{M} \phi_m(\tau,z) e\left( Nm\omega\right) \vert \le A_1$.  
Let $B_1$ be an open Euclidean ball centered at~$\Omega_o$ with closure 
$\overline{B}_1 \subseteq \Half_2$.  We can push the closed ball $\overline{B}_1$ down 
by $-i\smallmat000\epsilon$ so that this translation is still contained in~$\Half_2$.  
We do this in detail.  For any~$\epsilon>0$, the translated ball lies in the space 
$L=\{ \tzzw \in M_{2\times 2}^{\text{\rm sym}}(\C): \tau \in \Half_1\}$.  
Define the function $h:L\to\R$ for any $\Omega=\tzzw \in L$ by 
$$
h(\Omega) = 
\frac{\det(\IM \Omega)}{\langle \IM \Omega,\smallmat1000 \rangle }
= \IM \omega - \frac{(\IM z)^2}{\IM \tau}.  
$$
The function~$h$ is continuous on~$L$ and positive on~$\Half_2 \subseteq L$ 
and so has a positive minimum on~$\overline{B}_1$.  
Set $\epsilon=\frac12 \min_{\Omega \in \overline{B}_1} h(\Omega)>0$.  
We remark that~$\epsilon$ only depends on the choice of the ball~$B_1$.  
Translate the ball~$B_1$ by $-i\smallmat000\epsilon$ to obtain the ball 
$B_2=\{\Omega -i\smallmat000\epsilon \in L: \Omega \in B_1\}$ with closure in~$L$ 
given by  $\overline{B}_2=\{\Omega -i\smallmat000\epsilon \in L: \Omega \in \overline{B}_1\}$.  
We have 
$$
\inf_{\Omega \in \overline{B}_2}h(\Omega) = 
\inf_{\Omega \in \overline{B}_1}\left( h(\Omega) -\epsilon\right)  =
2\epsilon - \epsilon=\epsilon >0.  
$$
Therefore $\det(\IM \Omega)>0$ for $\Omega \in \overline{B}_2$ and 
$\overline{B}_2 \subseteq \Half_2$. By this process, $B_2$ is completely 
determined by~$B_1$.  
Next define a compact set 
$K_2=\{\Omega + \smallmat000\mu \in \Half_2: \Omega \in \overline{B}_2, 0 \le \mu \le 1 \}$ 
that contains $\overline{B}_2$.  
Use the continuity of the~$g_j$ to define 
$A_0=A_0(B_1,P)= \sup_{\Omega \in K_2} 
\left( 1+ \sum_{j=1}^d \vert g_j(\Omega)\vert \right)$.  
We will show that 
$A_1=\frac{A_0}{e^{2\pi N\epsilon}-1}$ %
works as the local bound 
at~$B_1$.  
\smallskip 

The main step will be to show 
$\vert \sum_{m=1}^{M} \phi_m(\tau_1,z_1) e\left( Nm\omega_1\right) \vert \le A_1$ 
for every $M \in \N$ and every 
$\Omega_1=\smallmat{\tau_1}{z_1}{z_1}{\omega_1} \in B_1 \cap \WW(\Gamma)$.  
The conclusion 
$\vert \sum_{m=1}^{M} \phi_m(\tau,z) e\left( Nm\omega\right) \vert \le A_1$ 
on the neighborhood $B_1$ 
then follows because 
$ \sum_{m=1}^{M} \phi_m(\tau,z) e\left( Nm\omega\right) $ is continuous 
on $ B_1$ and $B_1 \cap \WW(\Gamma)$ is dense in~$B_1$.  
The assumption $\Omega_1\in B_1 \cap \WW(\Gamma)$ implies $(\tau_1,z_1) \in \WW_1(\Gamma)$.  
We may rewrite the hypothesis 
$$
\fracg^d+
\cdots + \FJ(g_j) \fracg^{d-j}+ \cdots + \FJ(g_d)=0 
$$
as $\sum_{j=0}^d \FJ(g_j) \fracg^{d-j}=0$ where $g_0=1 \in M_0(K(N))^{+}$.  
Apply the ring homomorphism $\Atz(\tau_1,z_1): \FMpm(N) \to \C[[\xi]]$ to the given relation 
$\sum_{j=0}^d \FJ(g_j) \fracg^{d-j}=0$ in $\FM(kd,N,\epsilon^d;\Gamma)$ 
to obtain a relation among formal power series 
\begin{equation}
\label{eqatzg}
\sum_{j=0}^d \left(\Atz(\tau_1,z_1)\FJ(g_j)\right) \left(\Atz(\tau_1,z_1)\fracg\right)^{d-j}=0.  
\end{equation}
By Lemma~\ref{lemAtzs0}, we know $\Atz(\tau_1,z_1)\fracg \in \Sps(D_1)$ for % 
$\eta_1=\frac{(\IM {z_1})^2}{\IM {\tau_1}}$,  the radius 
$r_1=e^{-2\pi \eta_1} \le 1$, and the disk $D_1=\{z\in\C: \vert z \vert <r_1\}$. 
Moreover, the power series $\Atz(\tau_1,z_1)\FJ(g_j)$ converges in~$D_1$.    
Therefore the formal power series in equation~{(\ref{eqatzg})} is in the subring $\Sps(D_1)$  
of {\em convergent\/} power series on~$D_1$.  
\smallskip

For $\omega \in \Horo(\eta_1)$, we have $\xi_1=e(\omega) \in D_1$ and so we may 
use the evaluation homomorphism $[\xi_1]: \Sps(D_1) \to \C$ to obtain a relation among 
complex numbers
\begin{equation}
\label{eqatzgxi}
\sum_{j=0}^d \left(\Atz(\tau_1,z_1)\FJ(g_j)\right)[\xi_1] \left(\left(\Atz(\tau_1,z_1)\fracg\right)[\xi_1]\right)^{d-j}=0.  
\end{equation}
For $\omega \in \Horo(\eta_1)$, define 
$\Omega_1(\omega)=\smallmat{\tau_1}{z_1}{z_1}{\omega} \in \Half_2$.    
By Lemma~\ref{lemAtzs0},   for any $\omega \in \Horo(\eta_1)$, we have both 
$\left(\Atz(\tau_1,z_1)\FJ(g_j)\right)[\xi_1]=g_j\left( \Omega_1(\omega) \right)$ and 
$\left( \Atz(\tau_1,z_1) \fracg \right)[\xi_1]=\Hatz(\tau_1,z_1,\fracg)(\omega)$.  
We rewrite equation~{(\ref{eqatzgxi})} as 
\begin{equation*}
\label{eqatzgxiH}
\left(\Hatz(\tau_1,z_1,\fracg)(\omega)\right)^{d}  +
\sum_{j=1}^d g_j\left( \Omega_1(\omega) \right) \left(\Hatz(\tau_1,z_1,\fracg)(\omega)\right)^{d-j}=0.  
\end{equation*}

By Lemma~\ref{lempoly} we have 
\begin{equation}
\label{eqboundg}
\vert \Hatz(\tau_1,z_1,\fracg)(\omega)  \vert \le 
1+ \sum_{j=1}^d \vert g_j\left(\Omega_1(\omega) \right) \vert .  
\end{equation} 
By the definition of $A_0= \sup_{\Omega \in K_2} \left( 1+ \sum_{j=1}^d \vert g_j(\Omega)\vert \right)$, 
we have 
\begin{equation}
\label{eqboundgA0}
\vert \Hatz(\tau_1,z_1,\fracg)(\omega)  \vert \le A_0,  
\quad \text{if $\Omega_1(\omega)  \in K_2$.}
\end{equation} 
We have 
$\omega_1 -i \epsilon \in \Horo(\eta_1)=\{\omega \in \Half_1: \Omega_1(\omega) \in \Half_2\}$ because   
$\Omega_1\left( \omega_1 -i \epsilon \right) \in B_2 \subseteq \Half_2$.  
Furthermore, for $0\le \mu \le 1$, we have 
$\omega_1 -i \epsilon + \mu \in \Horo(\eta_1)$ because   
$\Omega_1\left( \omega_1 -i \epsilon+\mu \right) \in K_2 \subseteq \Half_2$, 
and  the function $\Hatz(\tau_1,z_1,\fracg)$ is holomorphic on 
$\omega_1 -i \epsilon+[0,1] \subseteq \Horo(\eta_1)$.  
Thus we may compute the Fourier coefficients of 
$$
\Hatz(\tau_1,z_1,\fracg)(\omega) =
\sum_{m=1}^{\infty} \phi_m(\tau_1,z_1) e\left( Nm\omega\right)
$$
by the familiar formula 
$$
\phi_m(\tau_1,z_1) = 
\int_0^1 \Hatz(\tau_1,z_1,\fracg)(\omega_1 -i \epsilon+\mu) 
e\left(-Nm(\omega_1 -i \epsilon+\mu)  \right)\,d\mu.
$$
We calculate, for any $M \in \N$, 
\begin{align*}
&\vert \sum_{m=1}^{M} \phi_m(\tau_1,z_1) e\left( Nm\omega_1\right) \vert =\\ 
&\vert \sum_{m=1}^{M} \left(\int_0^1 \Hatz(\tau_1,z_1,\fracg)(\omega_1 {-}i \epsilon{+}\mu)
e\left(-Nm(\omega_1 {-}i \epsilon{+}\mu)  \right)\,d\mu \right) e\left( Nm\omega_1\right) \vert \\
&\le 
\sum_{m=1}^{M} 
e^{-2\pi Nm\, \IM {\omega_1}} 
\int_0^1 \vert \Hatz(\tau_1,z_1,\fracg)(\omega_1 -i \epsilon+\mu) \vert 
e^{2\pi Nm (\IM {\omega_1} - \epsilon)}\,d\mu   \\
&= 
\sum_{m=1}^{M} 
e^{-2\pi Nm \epsilon} 
\int_0^1 \vert \Hatz(\tau_1,z_1,\fracg)(\omega_1 -i \epsilon+\mu)\vert 
\,d\mu.     
\end{align*}
As mentioned above, we have 
$\Omega_1\left( \omega_1 -i \epsilon+\mu \right) 
=\smallmat{\tau_1}{z_1}{z_1}{\omega_1 -i \epsilon+\mu}\in K_2$ for $0\le \mu\le 1$,  
so that 
$\vert \Hatz(\tau_1,z_1,\fracg)(\omega_1 -i \epsilon+\mu) \vert \le A_0$ 
by equation~{(\ref{eqboundgA0})}.  Thus 
\begin{align*}
\vert &\sum_{m=1}^{M} \phi_m(\tau_1,z_1) e\left( Nm\omega_1\right) \vert \\ 
\le A_0 &\sum_{m=1}^{M} 
e^{-2\pi Nm \epsilon}   
< \dfrac{A_0 e^{-2\pi N \epsilon}}{1-e^{-2\pi N \epsilon}} 
=    \dfrac{A_0}{e^{2\pi N \epsilon}-1} =A_1.  
\end{align*}
This completes the main step and, as explained, consequently shows that 
$\vert \sum_{m=1}^{M} \phi_m(\tau,z) e\left( Nm\omega\right) \vert \le A_1$ 
on the neighborhood $B_1$ of~$\Omega_o$ for all $M \in \N$.     
Since $\Omega_o \in \Half_2$ was arbitrary, 
 the sequence of partial sums 
$\sum_{m=1}^{M} \phi_m(\tau,z) e\left( Nm\omega\right)$ 
is locally bounded on~$\Half_2$.  
\end{proof}

%%%%%%%%%%%%%%%%%%%%%%%%%%%
%%%%%%%%%%%%%%%%%%%%%%%%%%%
\section{Holomorphicity.}
\label{secholomorphicity}
%%%%%%%%%%%%%%%%%%%%%%%%%%%
%%%%%%%%%%%%%%%%%%%%%%%%%%%

In Theorem~\ref{theoremord} we show that the divisibility of paramodular Fricke eigenforms 
is implied by the cuspidality of the quotient of their formal series.  Later, after proving our main result, 
Theorem~\ref{theoremordconclusion} improves Theorem~\ref{theoremord} by dropping the 
cuspidality assumption.

We review divisor theory following Gunning~\cite{gunning}.  
Let $G \subseteq \C^n$ be a domain.  
For open $U \subseteq G$, let $\OO(U)$ be the ring of holomorphic functions on~$U$.  
For $p \in G$, let $\OO_p$ be the ring of germs of holomorphic functions, and~$\MMM_p$ 
the field of germs of meromorphic functions at~$p$.  
A {\em holomorphic subvariety\/} is a subset $V \subseteq G$ such that each point $p\in V$ 
has a neighborhood~$U$ where $U \cap V$ is the zero set of finitely many functions holomorphic on~$U$.  
A holomorphic subvariety~$V$ is {\em irreducible\/} if $V = V_1 \cup V_2$ for holomorphic subvarieties $V_1$, $V_2$, 
implies $V_1=V$ or~$V_2=V$.  
The {\em regular points\/} of~$V$,  $\RR(V)$, are the points $p \in V$ where $U \cap V$ is a complex manifold 
for some neighborhood~$U$ of~$p$.  The regular points are an open dense subset of~$V$, 
locally finitely connected, and $\RR(V)$ is connected precisely when~$V$ is irreducible.   
The dimension, $\dim V$, of an irreducible~$V$ is the dimension of the complex manifold~$\RR(V)$, 
and the codimension of~$V$ is $n -  \dim V$.  
\smallskip

We let $\EE(G)$ be the set of irreducible, codimension one holomorphic subvarieties of~$G$.  
For a function $\Delta: \EE(G) \to \Z$, the {\em support\/} of~$\Delta$ is 
$\supp(\Delta)= \{ V \in \EE(G): \Delta(V) \ne 0 \}$.  
We say $\Delta$ is {\em locally finite\/} if, for all open $U \subseteq G$ with ${\bar U} \subseteq G$, 
the set $ \{ V \in \supp(\Delta): U \cap V \ne \emptyset \}$ is finite.  
The {\em group of divisors\/} in~$G$, $\DD_G$, 
is the abelian group of locally finite functions $\Delta: \EE(G) \to \Z$.  
Since $\supp(\Delta)$ is countable, 
we often write a global divisor as $\Delta = \sum_{j \in \N} \nu_j V_j$,  
for $\nu_j = \Delta(V_j) \in \Z$ and distinct $V_j \in \supp(\Delta)$.     
The semigroup of {\em effective\/} divisors, $\DD_G^{+} \subseteq \DD_G$, is defined by $\nu_j \ge 0$.   
\smallskip

For the germ of a holomorphic subvariety at~$p$, $V_p$, the dimension, $ \dim V_p$, 
is the maximal dimension of the finitely many connected components of~$\RR(V)$ whose closure contains~$p$.  
We say $V_p$ is {\em pure dimensionsal\/} when all these connected components have the same dimension.  
The definition of irreducibility for~$V_p$ is as before.  Let $\EE(p)$ be the set of germs of holomorphic subvarieties 
at~$p$ that are irreducible and have codimension one.  
A function $\delta: \EE(p) \to \Z$ has {\em support\/}  
$\supp(\delta)= \{ V_p \in \EE(p): \delta(V_p) \ne 0 \}$.  
The group of {\em local divisors\/} at~$p$, $\DD_p$, is the free abelian group on~$\EE(p)$.  
We write local divisors $\delta \in \DD_p$ as $\delta = \sum_{j =1}^{\ell} \nu_j V_j$ 
with $\nu_j = \delta(V_j) \in \Z$ and distinct $V_j \in \EE(p)$ but, fundamentally, 
a local divisor at~$p$ is a function $\delta: \EE(p) \to \Z$ with finite support.  
The semigroup of {\em effective\/} local divisors, $\DD_p^{+} \subseteq \DD_p$, is defined by $\nu_j \ge 0$.   
If $p \in \RR(V)$ then $V_p$ is irreducible.  
For $V \in \EE(G)$ and general~$p\in V$, $V_p = V_p' \cup V_p'' \cup \cdots$ decomposes into a finite union of distinct 
$V_p', V_p'', \ldots \in \EE(p)$.  
\smallskip

The ring $\OO_p$ is a unique factorization domain and noetherian.  
The fundamental correspondence of algebraic geometry holds between germs of holomorphic subvarieties 
and germs of holomorphic functions at~$p$.  
The ideal $\id(V_p) \subseteq \OO_p$ consists of the germs that vanish on~$V_p$ and, 
for an ideal $\fraca \subseteq \OO_p$, $\loc(\fraca)$ is the germ of the holomorphic subvariety 
defined by the simultaneous vanishing of the elements of~$\fraca$.  
The germ $V_p$ is irreducible precisely when $\id(V_p)$ is prime, and $V_p$ has pure 
codimension one precisely when $\id(V_p)$ is principal.  Therefore, 
for $V \in \EE(G)$ and $p \in V$, there is a neigborhood~$U$ of~$p$ and a $\varpi \in \OO(U)$, with $\varpi_p$ prime in~$\OO_p$,  
such that $V \cap U=\{ Z \in U: \varpi(Z)=0\}$.  We call this a local equation at~$p$ and refer to~$\varpi$ 
as a uniformizer.  
As a consequence of Cartan's Theorem, see Theorem~{F6} of~\cite{gunning},  
there is a neighborhood $U'' \subseteq U$ of~$p$ such that
$\varpi_q$ generates $\id V_q$ at all points $q \in U'' \cap V$.  
Cartan's Theorem asserts the existence of a neighborhood 
$U' \subseteq U$ of~$p$ and of functions $h_1, \ldots, h_{\ell} \in \OO(U')$ 
such that  $\id V_q=\OO_q(h_{1q}, \ldots, h_{\ell q} )$ at all points $q \in U' \cap V$.  
Since $\varpi_p$ divides each $h_{jp}$ there is a neighborhood $U''\subseteq U'$ 
of~$p$ where 
$\id V_q=\OO_q(h_{1q}, \ldots, h_{\ell q} ) \subseteq \OO_q \varpi_q$ at all points $q \in U'' \cap V$. 
Hence $\id V_q= \OO_q \varpi_q$ and $\varpi_q$ is prime for all $q \in U'' \cap V$. 
\smallskip

For a ring $R$, let $R^{\times}$ be the group of units and $R^*=R\,{\setminus} \{0\}$.  
The sequence of semigroups 
$ 0 {\to} \left( \OO_p^{\times}, \cdot\right) {\to} \left( \OO_p^*, \cdot\right) 
\overset{\divv_p}{\to} \left( \DD_p^{+},+\right) {\to} 0$  is exact, 
where  $\divv_p(f_p)=\sum_{j=1}^{\ell} \nu_j \loc(\OO_p f_j)$  
is the semigroup homomorphism 
determined by the factorization $f_p = \text{\rm(unit)} f_1^{\nu_1} \cdots f_{\ell}^{\nu_{\ell}}$ in $\OO_p$ 
into powers of nonassociate irreducibles 
$f_j$,  unique up to 
order and units.  The rule $\divv_p\left( f_p / g_p \right) = \divv_p(f_p) - \divv_p(g_p) $ defines 
a group homomorphism $\divv_p: \left( \MMM_p^* , \cdot\right)  \to \left( \DD_p,+\right)$.  
The local divisor $\divv_p(f_p)$ is effective precisely when 
$f_p \in \OO_p^* \subseteq \MMM_p^*$.  
To each global divisor $\Delta \in \DD_G$ and point $p \in G$, we associate a local divisor,  
$ \germp(\Delta) = \sum_{V \in \supp(\Delta):\hskip0.02in p \in V} \Delta(V) \left( V_p' + V_p'' + \cdots \right) $, 
where $V_p = V_p' \cup V_p'' \cup \cdots$ is the finite decomposition of~$V_p$ into irreducibles.  
For each $f \in \OO(G)^*$, the global divisor $\divv(f)$ is the unique $\Delta \in \DD_G^{+}$ 
such that 
$\germp(\Delta)= \divv_p(f_p)$ for all $p \in G$.  
 The rule $\divv\left( f / g \right) = \divv(f) - \divv(g) $ extends $\divv$ to  
a group homomorphism $\divv: \left( \MMM(G)^* , \cdot\right)  \to \left( \DD_G,+\right)$, 
where the field of meromorphic functions $\MMM(G)$ is the quotient field of~$\OO(G)$.  
The global assertion that $\divv(f)$ is effective precisely when 
$f \in \OO(G)^* \subseteq \MMM(G)^*$ 
follows from the corresponding local assertion.    
\smallskip

For $V \in \EE(G)$ and $f \in \OO(G)^*$, a direct way to compute 
$\nu=\ord_V(f)=\left( \divv(f) \right)(V)$ is to take a regular point $p \in \RR(V)$ and a 
local equation at~$p$, $V \cap U=\{ Z \in U: \varpi(Z)=0\}$, so that $V_p$ is irreducible 
and $\varpi_p$ is prime with $\OO_p \varpi_p=\id(V_p)$.  
The factor  $\varpi_p^{\nu}$ in the unique factorization of~$f_p$ in~$\OO_p$ defines $\nu \in \N_0$.  
This~$\nu$ is independent of the choice of the regular point~$p$ and of the uniformizer~$\varpi$.
\smallskip

We use 
$\cupp=\smallmat{-1}{-N}{0}{1} \in \GL(2,\Z)$ 
and 
$\Cupp = \smallmat{\cupp}{0}{0}{\cupp^*} \in \KN$ in the following lemma.   
The usefulness of the following lemma is that a point on the $\KNp$-orbit 
of a codimension one irreducible holomorphic subvariety of~$\Half_2$ 
can be found where 
the uniformizer~$\varpi$ is regular in~$\omega$ at that point.  This allows us to use the 
Weierstrauss preparation theorem on~$\varpi$ in the proof of~Theorem~\ref{theoremord}.  
In Lemma~\ref{lemregular}, the term ``regular" regettably has two distinct meanings.  
A holomorphic~$\varpi$ is regular in~$\omega$ at $(\tau_o,z_o,\omega_o)$ means $\varpi(\tau_o,z_o,\omega) \not\equiv 0$ 
in all neighborhoods of~$\omega_o$,  
whereas a point~$p$ of a holomorphic subvariety~$V$ is regular if 
$V \cap U$ is a complex manifold for some neighborhood~$U$ of~$p$.
\vskip0.1in

\begin{lemma}
\label{lemregular}
Let $k,N \in \N$, and $\epsilon \in \{ \pm 1\}$.  
Let $f \in M_k\left( \KN \right)^{\epsilon}$.  
Let $V \subseteq \Half_2$ be an irreducible component of~$\divv(f)$. 
For any $g \in \KNp$, $V_g=g\langle V \rangle$ is also an irreducible component of~$\divv(f)$, 
and we have $\ord_V f = \ord_{V_g} f$.  

There is a $g \in \{ I, \mu_N, \Cupp\} \subseteq \KNp$, a 
regular point $p \in V_g$, 
a neighborhood~$U \subseteq \Half_2$ of~$p$, and a function $\varpi \in {\mathcal O}(U)$ 
with an irreducible germ $\varpi_p \in \OO_p$, and 
with $\varpi\tzzw$ regular in~$\omega$ at~$p$, such that 
$V_g \cap U=\{ Z \in U: \varpi(Z)=0\}$. 
\end{lemma}
\begin{proof}
Each $g \in \KNp$ is biholomorphic on~$\Half_2$ so~$V$ is an irreducible holomorphic 
subvariety of codimension one if and only if~$V_g$ is.  
The automorphy of~$f$ shows that $V$ is supported in~$\divv(f)$ if and only if~$V_g$ is, 
although we remark that $V$ and $V_g$ might be equal.  
Since $g$ is biholomorphic, we have $\ord_V f = \ord_{g\langle V \rangle} f\circ g\inv$. 
The automorphy of~$f$ gives us $(f \circ g\inv)(\Omega)=j(g\inv,\Omega)^k f(\Omega)$ 
so that~$f$ and~$f\circ g\inv$ differ by a multiplicative unit and 
$\ord_{g\langle V \rangle} f\circ g\inv=\ord_{V_g} f$.  
\smallskip

Every regular point $p \in V$ has a local defining equation for~$V$ given by 
$V \cap U=\{ Z \in U: \varpi(Z)=0\}$ for some neighborhood~$U$ of~$p$, 
and some  uniformizer $\varpi \in {\mathcal O}(U)$
with an irreducible germ $\varpi_p$ represented by~$(\varpi,U)$.  
We may select~$U$ so that every point of~$U \cap V$ is also a 
regular point of~$V$ because~$V$ is a complex manifold at~$p$.  
By Cartan's Theorem, we may also assume that a single function element $(\varpi,U)$ 
satisfies $\id V_q= \OO_q \varpi_q$ for all $q \in U \cap V$.  
For the proof of the second part, we break the discussion into two mutually 
exclusive cases.  
{\em Case~$I$.  Some regular point $p \in V$ has $\varpi\tzzw$ regular in~$\omega$ or~$\tau$ at~$p$.\/}   
If $\varpi\smallmat{\tau_o}{z_o}{z_o}{\omega} \not\equiv 0$ in~$\omega$ for some regular 
$p=\smallmat{\tau_o}{z_o}{z_o}{\omega_o}\in V$ then we take~$g=I$ to satisfy our conclusion.  
Whereas if $\varpi\smallmat{\tau}{z_o}{z_o}{\omega_o} \not\equiv 0$ in~$\tau$ 
we take~$g=\mu_N$ and consider 
$q=\smallmat{\tau_1}{z_1}{z_1}{\omega_1}=\mu_N\langle p \rangle= 
\smallmat{\omega_o/N}{-z_o}{-z_o}{N\tau_o} \in V_g$.  
The point $q \in g\langle U\rangle$ is a regular point of~$V_g$  because $p$ is a regular point 
of~$V$ and~$g$ is biholomorphic.  
The local uniformizer at~$q \in V_g$ may be represented by $(\varpi \circ g\inv, g\langle U \rangle)$.  
In this case 
$$
(\varpi \circ g\inv)\smallmat{\tau_1}{z_1}{z_1}{\omega} = 
\varpi \smallmat{\omega/N}{-z_1}{-z_1}{N\tau_1}= 
\varpi \smallmat{\omega/N}{z_o}{z_o}{\omega_o}
$$
is regular in~$\omega$ at~$q$.  
From $V \cap U=\{ Z \in U: \varpi(Z)=0\}$ we deduce 
$V_g \cap g\langle U\rangle=\{ Z \in g\langle U\rangle: (\varpi\circ g\inv)(Z)=0\}$ so that~$q$ 
and $g\langle U\rangle$ play the role of~$p$ and~$U$ in the statement of the lemma.  
{\em Case II.  Every regular point $p \in V$ has 
$\varpi\tzzw$ not regular in~$\omega$ and not regular in~$\tau$ at~$p$.\/} 
Consider any regular point $p_o=\smallmat{\tau_o}{z_o}{z_o}{\omega_o}\in V$.  
Take a polydisk $\Delta=\Delta_{11} \times \Delta_{12} \times \Delta_{22} \subseteq U$ about~$p_o$ and 
use the function element $(\varpi, \Delta)$.   
We will ultimately show that Case~{II} is very special and that 
$V \cap \Delta=\{ Z \in \Delta: z-z_o=0\}$.  
We claim that $\varpi \smallmat{\tau}{z_o}{z_o}{\omega}=0$ for all 
$\smallmat{\tau}{z_o}{z_o}{\omega} \in \Delta$.  
If not there is some $r_o=\smallmat{\tau_1}{z_o}{z_o}{\omega_1} \in \Delta$ with $\varpi(r_o) \ne 0$.  
Since $\varpi$ is not regular in~$\omega$ at~$p_o$ we have 
$\varpi\smallmat{\tau_o}{z_o}{z_o}{\omega} \equiv 0$ for $\omega$ in some neighborhood of~$\omega_o$, 
hence necessarily for all $\omega \in \Delta_{22}$.  Therefore 
$\varpi\smallmat{\tau_o}{z_o}{z_o}{\omega_1} = 0$.  
The setting of Case~{II} applies equally well to 
$r=\smallmat{\tau_o}{z_o}{z_o}{\omega_1} \in V \cap \Delta$,  
which is still a regular point of~$V$, 
so that~$\varpi$ is not regular in~$\tau$ at $r$, and 
$\varpi\smallmat{\tau}{z_o}{z_o}{\omega_1} \equiv 0$ for $\tau$ in a neighborhood of~$\tau_o$, and 
hence for all $\tau \in \Delta_{11}$. 
Thus $\varpi(r_o)=\varpi\smallmat{\tau_1}{z_o}{z_o}{\omega_1} = 0$, contrary to our supposition.  
Thus $\varpi \smallmat{\tau}{z_o}{z_o}{\omega}=0$ for all 
$\smallmat{\tau}{z_o}{z_o}{\omega} \in \Delta$.  Therefore~$\varpi$ vanishes on 
$\{ Z \in \Delta: z-z_o=0\}$ and so the irreducible $z-z_o$ divides~$\varpi$ in~$\Delta$; 
however, the germ~$\varpi_{p_o}$ of~$\varpi$ at~$p_o$ is irreducible and so~$\varpi$ 
and~$z-z_o$ differ by a multiplicative unit in~$\Delta$, after possibly shrinking~$\Delta$ further.     
Therefore in Case~{II} we have the possible but special circumstance $V \cap \Delta=\{ Z \in \Delta: z-z_o=0\}$.  
Without loss of generality we may adjust the selection of~$\varpi$ by a unit and assume 
that $\varpi\tzzw=z-z_o$ in~$\Delta$.  
Now set~$g=\Cupp$ and consider 
$q=\smallmat{\tau_1}{z_1}{z_1}{\omega_1}=\Cupp\langle p_o \rangle= 
\smallmat{\tau_o+2Nz_o+N^2\omega_o}{z_o+N\omega_o}{z_o+N\omega_o}{\omega_o} \in V_g$.  
The local uniformizer at~$q \in V_g$ may be represented by $(\varpi \circ g\inv, g\langle \Delta \rangle)$.  
For this choice we have  
\begin{align*}
(\varpi \circ g\inv)\smallmat{\tau_1}{z_1}{z_1}{\omega} &= 
\varpi \smallmat{\tau_1-2Nz_1+N^2\omega}{z_1-N\omega}{z_1-N\omega}{\omega} \\
&= 
\varpi \smallmat{\tau_o+N^2(\omega-\omega_o)}{z_o+N(\omega_o-\omega)}{z_o+N(\omega_o-\omega)}{\omega}
=N(\omega_o-\omega),
\end{align*}
which is regular in~$\omega$ of order~$1$ at~$q$.  
Thus the conclusion of the lemma holds for~$q$ 
and $g\langle \Delta\rangle$ in the role of~$p$ and~$U$.  
\end{proof}

\begin{theorem}
\label{theoremord}
Let $k,k_0,N \in \N$ and $\epsilon, \epsilon_0 \in \{ \pm 1\}$.  
Let $\Gamma \subseteq \Gamma^0(N)$ be a subgroup of finite index. 
Let $\fracf \in \FMcusp(k,N,\epsilon;\Gamma)$. 
Let $f_0 \in M_{k_0}\left( \KN \right)^{\epsilon_0}$ be nontrivial and 
$G \in M_{k_0+k}\left( \KN \right)^{\epsilon_0\epsilon}$.   
Let $V \subseteq \Half_2$ be an irreducible component of~$\divv(f_0)$.  

If $\FJ(G)= \fracf\, \FJ(f_0)$ in $\FMcusp(k_0+k,N,\epsilon_0\epsilon;\Gamma)$ then 
$\ord_V f_0 \le \ord_V G$.  
The meromorphic~$G/f_0$ is holomorphic and $G/f_0 \in M_{k}\left( \KN \right)^{\epsilon}$.   
\end{theorem}
\begin{proof}
By Lemma~\ref{lemregular} we can find a point $p \in V_g$, for some $g \in \KNp$, 
where the local defining equation $V_g \cap U=\{ Z \in U: \varpi(Z)=0\}$ has 
the following nice properties.  The point~$p$ is a regular point of~$V_g$, 
the uniformizer~$\varpi$ is regular in~$\omega$ at~$p$, and the germ $\varpi_p$ 
is prime in~$\OO_p$.  
By Lemma~\ref{lemregular},  
in order to prove $\ord_{V} f_0 \le \ord_{V} G$ 
it suffices to prove $\ord_{V_g} f_0 \le \ord_{V_g} G$, 
so we may rename~$V_g$ as~$V$ without loss of generality.  
%\smallskip

In this paragraph we outline the remainder of the proof.  
We will find a neighborhood~$U_1$ of~$p$ such that every point $q \in U_1 \cap V$ 
has the same nice properties that~$p$ does.  We use this to twice move to nearby points.  
First we select $p_1 \in U_1 \cap V$ so that $\varpi_{p_1}$ is the only irreducible factor of $G_{p_1}$ 
that vanishes at~$p_1$.  Then we use the Weierstrauss preparation theorem 
on $\varpi_{p_1}$ to find a 
$p_2 \in \WW(\Gamma) \cap V \cap U_1$, 
and a holomorphic curve~$\Omega_2$ inside $\WW(\Gamma) \cap V $ and passing through~$p_2$,
which guarantees that the formal series~$\fracf$ converges on~$\Omega_2(\omega)$.  
The conclusion about the orders of vanishing readily follows from the convergence of~$\fracf$ 
on the holomorphic curve.    
\smallskip

The set of regular points $\RR(V)$ is open in~$V$.  
Regularity in~$\omega$ at~$p$ is a local condition in~$p$ 
because regularity of order less than or equal to~$\nu$ at~$p$ 
is implied by the nonvanishing of a partial derivative ${\partial^{\nu}_{\omega}\varpi}(p)$. 
By Cartan's Theorem, 
there is a neighborhood $U'' \subseteq U$ of~$p$ such that
$\varpi_q$ generates $\id V_q$ at all points $q \in U'' \cap V$.  
Accordingly there is a smaller neighborhood~$U_1$ of~$p$ where every $q \in U_1 \cap V$ 
has the nice properties $q \in \RR(V)$, $\varpi$ regular in~$\omega$ at~$q$, and $\varpi_q$ 
prime in~$\OO_q$.  
\smallskip  

In the local ring ${\mathcal O}_{p}$, factor 
$(f_0)_{p}= (\Delta_1)_{p} \varpi_{p}^{\nu(f_0)}$ and $G_{p}= (\Delta_2)_{p} \varpi_{p}^{\nu(G)}$,  
where $(\Delta_1)_{p},(\Delta_2)_{p}$ are 
finite products of irreducibles that are not associate to~$\varpi_{p}$, and where 
$\nu(f_0)=\ord_V f_0$ and $\nu(G) =\ord_V G$ by the direct way of computing the vanishing order on~$V$.  
We select a neighborhood $N_1 \subseteq U_1$ where the germs of the above factors all 
have representative function elements.  
If $\Delta_2(q)=0$ for all $q \in N_1 \cap V$ then $(\Delta_2)_p$ vanishes on~$V_p$ 
and~$\varpi_{p}$ divides $(\Delta_2)_p$ in~$\OO_p$ because~$\varpi_{p}$ is prime,  
contradicting the fact that $(\Delta_2)_p$ is a finite product of irreducibles that are not associate 
to~$\varpi_{p}$ in~$\OO_p$.  Hence there is a point $p_1 \in V \cap N_1$ such that $\Delta_2(p_1) \ne 0$.  
By shrinking $N_1$, we may assume that $\Delta_2$ is nonzero on~$N_1$.  
\smallskip

We apply the Weierstrauss preparation theorem to~$\varpi$, which is regular in~$\omega$ at 
$p_1=\smallmat{\tau_1}{z_1}{z_1}{\omega_1}$.  
There exist neighborhoods $U_2$ of $(\tau_1,z_1)$ and $W_2$ of $\omega_1$ such that, 
$N_2= U_2 \times W_2 \subseteq N_1$ and 
\begin{align*}
&\{Z =\tzzw \in N_2: \varpi(Z)=0\}
=  \\
&\{ Z =\tzzw \in N_2: \Wpoly(\tau, z, \omega):=\sum_{j=0}^t h_j(\tau,z) (\omega-\omega_1)^{t-j} =0\}, 
\end{align*}
for some $t \in \N$, $h_j \in  {\mathcal O}(U_2)$ for $j=0, \ldots,t$, and $h_0=1$.  
The set~$\WW_1(\Gamma)$ is dense in~$U_2$ by Lemma~\ref{lemdensity}, so by 
taking $(\tau_2,z_2)  \in \WW_1(\Gamma)$ close enough to $(\tau_1,z_1)\in U_2$ 
so that we may choose  
a root~$\omega_2$ of $\Wpoly(\tau_2,z_2,\omega)$ close enough to~$\omega_1$, 
we have 
$$
p_2=\smallmat{\tau_2}{z_2}{z_2}{\omega_2}   \in 
V \cap N_2 \cap \WW(\Gamma).  
$$
From $N_2 \subseteq N_1$ we inherit the representative function elements 
$f_0$, $G$, $\Delta_1$, $\Delta_2$, and~$\varpi$, so that 
$f_0= \Delta_1 \varpi^{\nu(f_0)}$ and $G= \Delta_2 \varpi^{\nu(G)}$ in $N_2$,  
and $\Delta_2$ is nonzero on~$N_2$. 

The formal series of~$G$ is given by 
$\FJ(G)= \fracf\, \FJ(f_0)$.  
Recalling that $(\tau_2,z_2)  \in \WW_1(\Gamma)$, 
and using Lemma~\ref{lemAtzs0} to apply $\Atz(\tau_2,z_2)$, we have the equality 
$\Atz(\tau_2,z_2)\FJ(G)= (\Atz(\tau_2,z_2)\fracf) (\Atz(\tau_2,z_2)\FJ(f_0))$ of  
convergent power series on~$D_2=\{\xi \in \C: \vert \xi \vert < r_2\}$ 
for $r_2=e^{-2\pi \eta_2}$ and $\eta_2=(\IM {z_2})^2/{\IM {\tau_2}}$.  
Also set $\Omega_2(\omega)=\smallmat{\tau_2}{z_2}{z_2}{\omega} \in \Half_2$ for $\omega \in \Horo(\eta_2)$.  
For $\omega \in \Horo(\eta_2)$ we have $e(\omega) \in D_2$ and may apply the evaluation homomorphism 
$[e(\omega)]: \Sps(D_2) \to \C$ to obtain 
$$
\left( \Atz(\tau_2,z_2)\FJ(G) \right)[e(\omega)]
= \left(\Atz(\tau_2,z_2)\fracf \right)[e(\omega)]\, 
\left(\Atz(\tau_2,z_2)\FJ(f_0) \right)[e(\omega)].  
$$
Making use of Lemma~\ref{lemAtzs0}, this may be written on 
$W_2 \cap \Horo(\eta_2)$, which contains a neighborhood of~$\omega_2$, as  
$$
G\left( \Omega_2(\omega) \right) = 
\Hatz(\tau_2,z_2,\fracf) (\omega)
f_0\left( \Omega_2(\omega) \right).  
$$
By specializing the factorizations of~$G$ and~$f_0$ in 
${\mathcal O}(N_2)$ to the holomorphic curve 
$\Omega_2(\omega)=\smallmat{\tau_2}{z_2}{z_2}{\omega} $ we have 
$$
\Delta_2 \left( \Omega_2(\omega) \right)
\varpi
\left( \Omega_2(\omega) \right)^{\nu(G)}= 
\Hatz(\tau_2,z_2,\fracf)(\omega) \,
\Delta_1 \left( \Omega_2(\omega) \right)
\varpi
\left( \Omega_2(\omega) \right)^{\nu(f_0)}.  
$$
Note that $\varpi\left( \Omega_2(\omega) \right)$ has at most~$t$ zeros on~$W_2$ 
and is hence nontrivial.  
Of the two cases $\nu(G) \ge \nu(f_0)$, and $\nu(G) < \nu(f_0)$, the first is our conclusion, so 
we will conclude the proof by showing that the second does not occur.  In the second case, by cancelling 
powers of the nontrivial $\varpi\left( \Omega_2(\omega) \right)$, 
we have 
$$
\Delta_2 \left( \Omega_2(\omega) \right)
= 
\Hatz(\tau_2,z_2,\fracf)(\omega) \,
\Delta_1 \left( \Omega_2(\omega) \right)
\varpi
\left( \Omega_2(\omega) \right)^{\nu(f_0)-\nu(G)}.  
$$
We evaluate these at $\omega_2\in W_2 \cap \Horo(\eta_2)$ to obtain 
$$
\Delta_2 (p_2)
= 
\Hatz(\tau_2,z_2,\fracf)(\omega_2) \,
\Delta_1 (p_2)
\left( 0^{\,\nu(f_0)-\nu(G)} \right) 
 =0,  
$$
which contradicts $\Delta_2 (p_2) \ne 0$. Thus we have $\nu(G) \ge \nu(f_0)$.  
\smallskip

A meromorphic function with an effective divisor is holomorphic.  
Therefore $G/f_0$ is holomorphic and continuous on~$\Half_2$.  
There is an open dense subset of~$\Half_2$ 
where $G/f_0$ transforms like an element of $M_{k}\left( \KN \right)^{\epsilon}$.  
By the continuity of $G/f_0$ and of the factor of automorphy, 
$G/f_0$ transforms like an element of $M_{k}\left( \KN \right)^{\epsilon}$ on~$\Half_2$ 
and hence  $G/f_0 \in M_{k}\left( \KN \right)^{\epsilon}$.  
\end{proof}

%%%%%%%%%%%%%%%%%%%%%%%%%%%
%%%%%%%%%%%%%%%%%%%%%%%%%%%
\section{Main Theorem.}
\label{secmain}
%%%%%%%%%%%%%%%%%%%%%%%%%%%
%%%%%%%%%%%%%%%%%%%%%%%%%%% 

\begin{theorem}
\label{theoremSVOdense} 
Let $U \subseteq \C^d$ be open.  
Let $\{ f_j\}$ be a locally bounded sequence of holomorphic functions 
on~$U$ that converges on a dense subset of~$U$.  
Then the sequence $\{ f_j\}$ converges on~$U$ and uniformly 
on compact subsets of~$U$.   
\end{theorem}
\begin{proof}
This is Exercise $4$a of section~$4$ in Chapter~$1$ of~\cite{MR1893803}.  

The proof given here imitates the proof of Lemma~{IV.4.8} in~\cite{MR2513384} 
for the one-dimensional case.  
It suffices to show that $\{ f_j\}$ is uniformly Cauchy on compact subsets of~$\Omega$.   
Pick a compact $K \subseteq \Omega$ and an $\epsilon>0$.  
There is no loss of generality in assuming that~$K$ is a closed ball.  
We will construct an $M \in \N$ from the given data $\{ f_j\}$, $K$, and $\epsilon$.  

The family $\{ f_j\}$ is equicontinuous on~$K$ because it is locally bounded, 
see pages~$11$-$12$ of~\cite{MR3526579}.    
Select $\delta>0$ to enforce this equicontinuity on~$K$ for~$\frac13 \epsilon$.  
\begin{equation}
\label{eqequicont}
\forall z,w\in K,\, \forall n \in \N,\, 
\vert z-w \vert < \delta \implies \vert f_n(z) - f_n(w) \vert < \tfrac13 \epsilon.  
\end{equation}
Cover~$K$ with open balls inside~$\Omega$ at each point of~$K$ with radii that are less than~$\frac12 \delta$.  
By the compactness of~$K$ we have 
$K \subseteq \cup_{i=1}^{\ell} B(z_i, \delta_i) \subseteq \Omega$ for some 
$z_i \in K$, and $0< \delta_i < \frac12 \delta$.  
Since~$S$ is dense in~$K$, we can pick an $s_i \in K \cap S \cap B(z_i, \delta_i)$; 
this uses that~$K$ is a closed ball.   
Note that for every point $z \in B(z_i, \delta_i)$ we have $ \vert z-s_i \vert < \delta$ and 
therefore $\vert f_n(z) - f_n(s_i) \vert < \frac13\epsilon$ by~{(\ref{eqequicont})}.  
By hypothesis each sequence~$f_j(s_i)$ converges and so is Cauchy.  
Select $M_i>0$ so that 
\begin{equation}
\label{eqcauchy}
\forall m,n \in \N,\, 
m,n>M_i \implies \vert f_m(s_i) - f_n(s_i) \vert <\tfrac13 \epsilon.  
\end{equation}
Let $M=\max_{1\le i \le \ell} M_i$ be the promised natural number.  
We have 
$$
\vert f_m(z) {-} f_n(z) \vert \le 
\vert f_m(z) - f_m(s_i) \vert + 
\vert f_m(s_i) - f_n(s_i) \vert+
\vert f_n(s_i) - f_n(z) \vert
$$
for any $m,n\in\N$, $z \in K$, and $1\le i \le \ell$.  
If $m,n>M$, then $m,n>M_i$ and so $\vert f_m(s_i) - f_n(s_i) \vert < \frac13 \epsilon$ by~{(\ref{eqcauchy})}.  
Choose~$i$ so that $z \in B(z_i, \delta_i)$; then 
$\vert f_j(z) - f_j(s_i) \vert < \frac13 \epsilon$ for all $j \in \N$ by~{(\ref{eqequicont})}.  
Therefore, for all $m,n\in\N$ with $m,n>M$, and all $z \in K$, we have 
$\vert f_m(z) - f_n(z) \vert<\epsilon$, and the sequence $\{ f_j\}$ is uniformly Cauchy on~$K$.  
\end{proof}
%\vskip0.2in

The essence of the proof of the main result, Theorem~\ref{thmain3},  
lies in the argument for the following special case of Fricke plus cusp forms. 
The general case will be reduced to Theorem~\ref{thmain1}.  

\begin{theorem}
\label{thmain1}
The map 
$\FJ: S_k\left( \KNp \right) \to \FMcusp(k,N,+)$ 
is an isomorphism for $N,k\in\N$.  
\end{theorem}
\begin{proof}
Since the map~$\FJ$ is injective, the only issue is surjectivity.  
Take a nontrivial $\fracf \in \FMcusp(k,N,+)$.  
By Proposition~\ref{propfiniteindex}, 
there is a subgroup $\Gamma$ of finite index in~$\Gamma^0(N)$ such that 
$\fracf \in \FM(k,N,+\,; \Gamma)$, which implies  
$\fracf \in \FM(k,N,+\,; \Gamma) \cap \FMcusp(k,N,+)=\FMcusp(k,N,+\,; \Gamma)$.  
Completely separately, 
by Proposition~\ref{propalgebraic},
$\fracf$ satisfies a polynomial relation in $\FM(k_0+dk,N,+)$,  
$$
\FJ(f_0) \fracf^d + \FJ(f_1)\fracf^{d-1} + \cdots +\FJ(f_j)\fracf^{d-j} + \cdots + \FJ(f_d)=0,
$$
for some $d \in \N$, $k_0\in \N_0$, and some 
$f_j \in M_{k_0+jk}(\KN^{+})$ with $f_0$ not identically zero.  
Since $\FJ(f_j)$, $\fracf \in \FM(N,+\,;\Gamma)$, this polynomial relation is actually in $\FM(k_0+dk,N,+\,;\Gamma)$.  
Since $\FMcusp(N,+\,;\Gamma)$ is a graded ideal in $\FM(N,+\,;\Gamma)$, set 
$$\fracg=\FJ(f_0) \fracf = \sum_{m=1}^{\infty} \phi_m \xi^{Nm} \in \FMcusp(k+k_0,N,+\,;\Gamma).$$    
Then $\fracg$ is integral over $\FJ\left( M(\KNp) \right)$ and
satisfies a monic polynomial relation  in $\FM((k+k_0)d,N,+\,;\Gamma)$ given by  
$$
\fracg^d+\FJ(g_1) \fracg^{d-1}+
\cdots + \FJ(g_j) \fracg^{d-j}+ \cdots + \FJ(g_d)=0 
$$
with $g_j=f_0^{j-1}f_j \in M_{(k+k_0)j}\left( \KNp \right)$ for $j=1,\ldots,d$.  

By Theorem~\ref{lemlocallybounded},  the sequence of partial sums 
$\sum_{m=1}^{M} \phi_m(\tau,z) e\left( Nm\omega\right)$ 
is locally bounded for $\tzzw\in \Half_2$.  By Corollary~\ref{coraokinew}, 
we also know that the formal series~$\fracg$ converges on the dense 
subset~$\WW(\Gamma)$ of~$\Half_2$, noting 
$\fracg \in \FMcusp(k+k_0,N,+\,;\Gamma) \subseteq \FMcusp(k+k_0,N;\Gamma)$. 
Therefore, by Theorem~\ref{theoremSVOdense}, 
the sequence of partial sums converges uniformly on compact sets of~$\Half_2$.  
The limiting function $G \in {\mathcal O}(\Half_2)$ is given by 
$G(\Omega)=\sum_{m=1}^{\infty} \phi_m(\tau,z) e\left( Nm\omega\right)$ 
and satisfies $G\vert_{k+k_0}g=G$ for all $g \in \PtwooneZ$ 
because each term $\phi_m(\tau,z) e\left( Nm\omega\right)$ is so invariant.  
The holomorphic function~$G$ is 
hence periodic with respect to the translation lattice 
%% $t\left( M_{2 \times 2}^{\text{\rm sym}}(\Z)\right)$   
$\{ \smallmat{I}{s}{0}{I} \in \PtwooneZ: s \in M_{2 \times 2}^{\text{\rm sym}}(\Z) \}$ 
and has a Fourier series that converges absolutley and uniformly on compact sets of~$\Half_2$.   
This absolute convergence of the Fourier series of~$G$ is the crux of the matter because 
it will allow us to rearrange the order of summation 
and thereby deduce the crucial invariance property $G\vert_{k+k_0} \mu_N=G$ 
from the involution condition $c(n,r;\phi_m)=(-1)^k\epsilon\, c(m,r;\phi_n)$ on~$\fracg$. 
Since the convergence of the sequence of partial sums is uniform 
on compact subsets, the Fourier coefficients of~$G$ are 
supported on $\XtwosN \subseteq M_{2 \times 2}^{\text{\rm sym}}(\Q)$ and, 
for $t =\nrNm\in \XtwosN$, satisfy $a(t;G)=c(n,r;\phi_m)$.   
For~$\fracg$ in this proof we have $\epsilon=+1$.  
Recall that $t \mapsto  F_N\inv t F_N^{*}$ is an involution of~$\XtwosN$ 
satisfying $F_N\inv\smallmat{n}{r/2}{r/2}{Nm}F_N^{*}=\smallmat{m}{-r/2}{-r/2}{Nn}$, and that 
absolutely convergent 
summations over $t \in \XtwosN$ may instead be rearranged to be taken over $F_N\inv t F_N^{*}$.   
\begin{align*}
&\left( G \vert \mu_N \right)(\Omega) = G(F_N^{*} \Omega F_N\inv) 
= \sum_{t \in \XtwosN} a(t;G) 
e\left(\langle F_N^{*} \Omega F_N\inv, t \rangle \right) \\
= &\sum_{t} a(t;G) 
e\left(\langle  \Omega, F_N\inv t F_N^{*}   \rangle \right)
= \sum_{t} a(F_N\inv t F_N^{*};G) 
e\left(\langle  \Omega, t   \rangle \right) \\
= &\sum_{t={\tiny \nrNm} \in \XtwosN} c(m,-r;\phi_n) 
e\left(\langle  \Omega, t   \rangle \right)
=
\sum_{t } c(n,r;\phi_m) 
e\left(\langle  \Omega, t   \rangle \right) \\
= &\sum_{t} a(t;G) 
e\left(\langle  \Omega , t \rangle \right) =G(\Omega).  
\end{align*}
Combining these automorphic properties of~$G$ we obtain 
$G\vert \sigma = G$ for all 
$ \sigma \in \langle \PtwooneZ,\mu_N \rangle=\langle \KN,\mu_N \rangle=\KNp$, 
see Gritsenko~\cite{Gritsenko1995,MR1336601}.      
Therefore $G \in M_{k+k_0}\left( \KNp \right)$ and $\FJ(G)=\fracg \in \FMcusp(k+k_0,N,+\,;\Gamma)$.  
Together these imply that $G \in S_{k+k_0}\left( \KNp \right)$ by    
Reefschl{\"a}ger's decomposition~\cite{reefschlager73}, as in Lemma~\ref{lemringone}.  
\smallskip

We have $\FJ(G)=\fracg= \FJ(f_0) \fracf$ 
for the formal series $\fracf \in  \FMcusp(k,N,+\,;\Gamma)$, 
and the paramodular forms 
$f_0  \in M_{k_0}\left( \KNp \right)$, 
 $G \in S_{k+k_0}\left( \KNp \right)$.  
By Theorem~\ref{theoremord}, 
$f{:=}\, {G}/{f_0} \in M_{k}\left( \KNp \right)$. 
From $f_0f=G$ we see that $\FJ(f_0) \FJ(f) = \FJ(G)$.  
Combining this with $\FJ(f_0) \fracf=\FJ(G)$, we have 
$\FJ(f_0)\left( \FJ(f)-\fracf\right)=0$, which implies $\FJ(f)=\fracf$ since 
$\FM(N)$ is an integral domain and $f_0$ is nontrivial.  
By Reefschl{\"a}ger's result, $\FJ(f)=\fracf \in \FMcusp(k,N,+\,;\Gamma)$ and $f \in M_{k}\left( \KNp \right)$ imply 
$f \in S_{k}\left( \KNp \right)$, which proves that $\FJ$ is surjective.   
\end{proof}

\begin{lemma}
\label{lemmacuspminus}
There is a nontrivial cusp form $\chi_{N} \in S_{11}\left( \KN \right)^{-}$ for $N >1$.     
\end{lemma}
\begin{proof}
There are a number of possible approaches to proving the existence of 
nontrivial minus cusp forms. We might use an asymptotic trace formula or 
the oldform theory of Roberts and Schmidt~\cite{robertsschmidt06,robertsschmidt07}.  
Perhaps the briefest is to use the injectivity of the Gritsenko lift~\cite{MR1336601}, 
$\Grit: \Jkmc \to S_{k}\left( K(m) \right)^{\epsilon}$ for $\epsilon=(-1)^k$.   %
An estimate from the 
dimension formula~\cite{ez85} is 
$\dim J_{11,m}^{\text{\rm cusp}} \ge \frac74 \floor\left(\frac{m-3}{6} \right) +\frac{7}{12}$ for $m \in \N$.    
Considering $\dim J_{11,2}^{\text{\rm cusp}}=1$,  
%%  and $\dim J_{11,3}^{\text{\rm cusp}}=1$, 
we see $\dim J_{11,m}^{\text{\rm cusp}} \ge 1$ for $m \ge 2$.  
Thus there is a nontrivial Gritsenko lift~$\chi_{N}$ in $S_{11}\left( \KN \right)^{-}$ for $N >1$. 
\end{proof}

\begin{theorem}
\label{thmain2}
The map 
$\FJ: S_k\left( \KN \right)^{\epsilon} \to \FMcusp(k,N,\epsilon)$ 
is an isomorphism for $N,k\in\N$, and $\epsilon \in \{ \pm 1\}$.    
\end{theorem}
\begin{proof}
The case $\epsilon=+1$ is Theorem~\ref{thmain1}, so assume $\epsilon=-1$ and 
$N>1$.    
Take $\fracf \in \FMcusp(k,N,-)$.    
There is a nontrivial cusp form $\chi_{N} \in S_{11}\left( \KN \right)^{-}$ by 
Lemma~\ref{lemmacuspminus}.  Multiply $\fracf$ by $\FJ(\chi_N)$ to obtain a new  
formal series $\fracg=\FJ(\chi_N) \fracf \in \FMcusp(k+11,N,+)$.  
By Theorem~\ref{thmain1}, there exists a 
$G \in S_{k+11}\left( \KN \right)^{+}$ such that $\FJ(G)=\fracg=\FJ(\chi_N)\fracf$.  
Both $\FJ(G)$ and $\FJ(\chi_N)$ have $\Gamma^0(N)$-symmetries and we will show that 
$\fracf$ does too.  
Apply the ring homomorphism $\AFS: \FM(N) \to \ffs$ of section~\ref{secinvariance} 
to $\FJ(G)=\FJ(\chi_N) \fracf$ to obtain 
$\FS(G)=\FS(\chi_N) \AFS(\fracf)$.  
Now, for $\sigma \in \Gamma^0(N)$, apply the automorphism~$\jj(\sigma)$ 
of Lemma~\ref{lemmajj},  
which fixes $\FS(G)$ and $\FS(\chi_N)$ to obtain 
$0=\FS(\chi_N) \left(\AFS(\fracf)-\jj(\sigma)\AFS(\fracf) \right)$.  
Since $\ffs$ is an integral domain and $\FS(\chi_N)$ is nontrivial, 
we have $\jj(\sigma)\AFS(\fracf)=\AFS(\fracf)$ for all $\sigma \in \Gamma^0(N)$.  
By Lemma~\ref{lemmajj}, we conclude that $\fracf \in \FMcusp(k,N,-\,;\Gamma^0(N))$.  %%  
By Theorem~\ref{theoremord} with $\Gamma=\Gamma^0(N)$, 
the equation $\FJ(G)=\FJ(\chi_N) \fracf$ implies that $f{:=}\,G/\chi_N$ is holomorphic, 
and, consequently, 
that $f \in M_k\left( \KN \right)^{-}$.  
From $\FJ(G)=\FJ(\chi_N) \fracf$ and $\FJ(G)=\FJ(\chi_N) \FJ(f)$ we obtain 
$0=\FJ(\chi_N)(\fracf-\FJ(f))$. 
Since $\FM(N)$ is an integral domain and $\FJ(\chi_N)$ is nontrivial, we have $\FJ(f)=\fracf$.  
By Reefschl{\"a}ger's result and $\fracf \in \FMcusp(k,N,-\,;\Gamma)$ 
we have $f \in S_k\left( \KN \right)^{-}$.  
Thus $\FJ$ is surjective.  
\end{proof}

\begin{lemma}
\label{lemmacuspnonvanish}
Let $N \in \N$ and $I=[\KNp:\KNp \cap K(1)]$.  
Let $p \in \Half_2$.  
There is a cusp form $\chi_p \in S_{60I}\left( \KNp \right)$ such that $\chi_p(p) \ne 0$.  
\end{lemma}
\begin{proof}
The divisor of $\chi_{10} \in S_{10}(K(1))$ is the reducible locus~$R$,
$$
R=\cup_{\sigma \in K(1)} 
\{ \sigma \langle\smallmat{\tau}00{\omega}\rangle \in\Half_2: \tau, \omega \in \Half_1 \}.
$$ 
For $\psi_{12} \in S_{12}(K(1))$, we have 
$\psi_{12}\smallmat{\tau}00{\omega}=\Delta(\tau)\Delta(\omega)$ for 
the nonvanishing 
$\Delta \in S_{12}\left( \SL(2,\Z) \right)$.  
Thus $\chi_{10}$ and~$\psi_{12}$ have no common zeros on~$\Half_2$.  
Pick $I$ coset representatives~$\sigma$ for 
$\KNp = \amalg_{\sigma} (\KNp \cap K(1)) \sigma$.  
For any $(\alpha,\beta) \in \C^2 \setminus \{(0,0)\}$, consider the nontrivial norm 
$$
H^{\alpha,\beta} = \prod_{[\sigma] \in (\KNp \cap K(1))\backslash \KNp} 
\left(\alpha \chi_{10}^6+\beta \psi_{12}^5\right)\vert_{60} \sigma \in  S_{60I}\left( \KNp \right).  
$$
We have $H^{\alpha,\beta}(p)= 
\left( \prod_{\sigma} j(\sigma,p)^{-60}  \right) 
\prod_{\sigma}
\left(\alpha \chi_{10}^6(\sigma\langle p \rangle)+\beta \psi_{12}^5(\sigma\langle p \rangle)\right)$.  
Pick any $(\alpha_o,\beta_o) \in \C^2$ that is not on any of the finitely many 
codimension one  
complex lines through the origin $\{ (z_1,z_2) \in \C^2: 
z_1 \chi_{10}^6(\sigma\langle p \rangle)+z_2 \psi_{12}^5(\sigma\langle p \rangle)=0\}$.  
The cusp form $\chi_p=H^{\alpha_o,\beta_o}$ is as claimed.  
\end{proof}

\begin{theorem}[Main Theorem] 
\label{thmain3}
Let $N\in\N$, $k \in \N_0$, and $\epsilon \in \{ \pm 1\}$. 
The map 
$\FJ: M_k\left( \KN \right)^{\epsilon} \to \FM(k,N,\epsilon)$ 
is an isomorphism.     
The map 
$\FJ: S_k\left( \KN \right)^{\epsilon} \to \FMcusp(k,N,\epsilon)$ 
is also an isomorphism. 
\end{theorem}
\begin{proof}
The case $k=0$ is easy so assume $k \in \N$.  
Take $\fracf \in \FM(k,N,\epsilon)$. 
For any point $p \in \Half_2$ there is a cusp form 
$\chi_p \in S_{60I}\left( \KNp \right)$ such that $\chi_p(p) \ne 0$, 
as in Lemma~\ref{lemmacuspnonvanish}.  
Define a corrsponding product of formal series $\fracg_p=\FJ(\chi_p) \fracf \in \FMcusp(k+60I,N,\epsilon)$.  
By Theorem~\ref{thmain2}, there exists a 
$G_p \in S_{k+60I}(\KN)^\epsilon$ such that $\FJ(G_p)=\fracg_p=\FJ(\chi_p) \fracf$.  
Note that $S_{k+60I}(\KN)^\epsilon= S_{k+60I}(\KNp, \lambda_{\epsilon})$ 
where $\lambda_{\epsilon}: \KNp \to \{ \pm 1\}$ is the character defined by 
$\lambda_{\epsilon}(\mu_N)=\epsilon$ and $\lambda_{\epsilon}=1$ on~$\KN$.  
\smallskip

In the same way, 
for any point $q \in \Half_2$ distinct from~$p$ we also have 
$\FJ(G_q)=\FJ(\chi_q) \fracf$. Since $\FM(N)$ is an integral domain we necessarily have 
$\FJ(G_p) \FJ(\chi_q) - \FJ(G_q) \FJ(\chi_p)=0$, or   
$G_p\chi_q = G_q\chi_p$. 
Therefore we have 
$\frac{G_p}{\chi_p}=\frac{G_q}{\chi_q}$ whenever both quotients are defined.  
We may define $f \in {\mathcal O}(\Half_2)$ by setting $f=\frac{G_p}{\chi_p}$ in neighborhoods 
where the denominator does not vanish, and by noting this is independent of the choice of~$p$.  
Therefore we have $\chi_p f = G_p$, $\chi_p\vert \sigma\, f\vert \sigma = G_p\vert \sigma$, and 
$\chi_p f\vert \sigma = \lambda_{\epsilon}(\sigma) G_p$, for any point~$p$ and any $\sigma \in \KNp$.  
Thus $\chi_p \left( f\vert \sigma - \lambda_{\epsilon}(\sigma) f\right)=0$.  
Since $M(\KN)$ is an integral domain and $\chi_p$ is nontrivial, we have $f\vert \sigma = \lambda_{\epsilon}(\sigma) f$ and 
$f \in M_{k}(\KNp, \lambda_{\epsilon})=M_{k}(\KN)^\epsilon$.  
It remains to show $\FJ(f)=\fracf$.  
This follows from $\FJ(G_p)=\FJ(\chi_p) \fracf$, $\FJ(G_p)=\FJ(\chi_p) \FJ(f)$, and the fact that 
$\FM(N)$ is an integral domain and $\chi_p$ is nontrivial.    
The result for cusp forms follows from the general result because, by Lemma~\ref{lemringtwo}, 
the inverse image of $\FMcusp(k,N,\epsilon)$ under $\FJ$ is $S_k\left( \KN \right)^{\epsilon}$.  
\end{proof}

To help prove Corollary~\ref{corintro} we use a lemma 
from linear algebra.  
\begin{lemma}
\label{lemmaLA}
Let $U$, $V$, and~$W$ be finite dimensional vector spaces.  
Let $L':U \to W$ and $L'':V\to W$ be linear maps.  
Define a subspace $X \subseteq U \times V$ by 
$X=\{ (u,v) \in U \times V: L'u=L''v\}$.  
Define the projection $\pi_U: X \to U$ by $\pi_U(u,v)=u$.  
If $L''$ is injective then $\pi_U$ is injective, and conversely.  
\end{lemma}
\begin{proof}
We have 
$\ker( \pi_U) =\{ (0,v) \in U \times V: 0=L''v\}= \{0\} \times \ker(L'')$.  
Therefore $\ker( \pi_U)=\{(0,0)\}$ if and only if $\ker(L'')=\{0\}$.  
\end{proof}
We now prove Corollary~\ref{corintro} from the Introduction.  
\smallskip

\noindent
{\em Proof of Corollary~\ref{corintro}.}  
To prove that the sequence $\dim \FM(k,N,\epsilon)[d]$ is monotonically decreasing for $d \ge \frac16 Nk$, 
it suffices to prove that the projection 
$\pi_d: \FM(k,N,\epsilon)[d+1] \to \FM(k,N,\epsilon)[d]$, 
sending $\left( \phi_m \right)_{m=0}^{d+1}$ to $\left( \phi_m \right)_{m=0}^{d}$, injects.  
To employ the notation of Lemma~\ref{lemmaLA}, we set $U=\FM(k,N,\epsilon)[d]$ 
and $V = J_{k, N(d+1)}$.  To define~$W$, first define the finite set 
$I=\{ (m,r) \in \Z^2: 0 \le m \le d,\, 4N(d+1)m \ge r^2 \}$ and set $W= \C^{\vert I \vert}$.  
Define a linear map $L':U \to W$ by 
$L'\left(\left( \phi_m \right)_{m=0}^{d} \right) = 
\left( c(d+1,r;\phi_m) \right)_{(m,r) \in I}$, and   a linear map 
$L'':V \to W$ by 
$L''\left( \phi_{d+1}  \right) = 
\left( (-1)^k \epsilon\,c(m,r;\phi_{d+1}) \right)_{(m,r) \in I}$.  
From the definition of~$L''$ we see that the kernel of~$L''$ consists of the 
Jacobi forms vanishing to order~$d+1$, so that 
$\ker(L'')=J_{k, N(d+1)}(d+1)$.  
From equation~{(\ref{eqinvcondintro})} we see that $\FM(k,N,\epsilon)[d] \subseteq X$ for 
$$
X=\{ (\left( \phi_m \right)_{m=0}^{d}  , \phi_{d+1} ) \in U \times V: L'\left(\left( \phi_m \right)_{m=0}^{d} \right)=L''(\phi_{d+1})\}.  
$$
Since $d+1 > \frac16 Nk$, Proposition~\ref{propSQuaREvanishing} implies 
$\ker(L'')=J_{k, N(d+1)}(d+1)=\{0\}$, so that $\pi_U: X \to U$ injects by Lemma~\ref{lemmaLA}.  
Since $\pi_d$ is the restriction of~$\pi_U$ to $\FM(k,N,\epsilon)[d+1]$, 
we have shown that the projection 
$\pi_d: \FM(k,N,\epsilon)[d+1] \to \FM(k,N,\epsilon)[d]$ injects for $d \ge \frac16 Nk$.   
\smallskip

% Don't like \varprojlim_{d}
For the inverse system $\pi_d: \FM(k,N,\epsilon)[d+1] \to \FM(k,N,\epsilon)[d]$ we have 
$\FM(k,N,\epsilon)= \varprojlim_{d} \FM(k,N,\epsilon)[d]$.  
The reason for this is that any equation $c(n,r;\phi_m) = (-1)^k \epsilon\, c(m,r;\phi_m)$ 
imposed by equation~{(\ref{eqinvcondintro})} on 
$\fracf = \left( \phi_m \right)_{m=0}^{\infty} \in \FM(k,N,\epsilon)$ is also imposed 
on $\pi_d\fracf = \left( \phi_m \right)_{m=0}^{d} \in \FM(k,N,\epsilon)[d]$ for 
 all $d \ge \max(m,n)$. 
\smallskip

Since the sequence $\dim \FM(k,N,\epsilon)[d] \in \N_0$ is eventually monotonically decreasing, 
it is eventually constant, and the injective~$\pi_d$ are eventually isomorphisms.  Thus 
$\dim \FM(k,N,\epsilon)= \lim_{d} \dim\FM(k,N,\epsilon)[d]$.  
From Theorem~\ref{thmain3} we have 
$\dim M_k\left( \KN \right)^{\epsilon}= \dim \FM(k,N,\epsilon)$ and this completes the proof.

%%%%%%%%%%%%%%%%%%%%%%%%%%% 
%%%%%%%%%%%%%%%%%%%%%%%%%%%
\section{Concluding Remarks.}
\label{secconclusions}
%%%%%%%%%%%%%%%%%%%%%%%%%%%
%%%%%%%%%%%%%%%%%%%%%%%%%%%

The Main Theorem~\ref{thmain3} gives an quick proof,  
without needing to check convergence, 
of the existence of the 
Gritsenko lift~\cite{MR1336601}, $\Grit: J_{k,N} \to M_k\left( \KN \right)^{\epsilon}$, for $\epsilon = (-1)^k$.  
For $\ell \in \N$, let $V_{\ell}:J_{k,m} \to J_{k,m\ell}$ be the index raising operator from~\cite{ez85}.  
For $\phi \in J_{k,N}$, the formal series 
$\fracf=c(0,0;\phi) \frac12 \zeta(1-k) E_k + 
\sum_{m \in \N} \left( \phi | V_m\right) \,\xi^{Nm}$ is directly checked to be in $\FM(k,N,\epsilon)$.  
Here, the first term with the elliptic Eisenstein series~$E_k$ is zero unless $k \ge 4$ is even. 
By the Main Theorem there is a $\Grit(\phi) \in M_k\left( \KN \right)^{\epsilon}$ with 
$\FJ\left(\Grit(\phi) \right)=\fracf$.  
\smallskip

The Main Theorem~\ref{thmain3} can be used to construct the global level raising operators in the 
paramodular newform theory of Roberts and Schmidt~\cite{robertsschmidt06,robertsschmidt07}.  
For a prime~$\ell$, the three level raising operators 
$\theta_{\ell},  \theta'_{\ell}: S_k\left( K(N)  \right) \to  S_k\left( K(N\ell)  \right)$, and 
$\eta_{\ell}: S_k\left( K(N)  \right) \to  S_k\left( K(N\ell^2)  \right)$, 
are used to create the oldforms.  We can directly define the global level raising operators 
by giving their action on Fourier-Jacobi expansions.  
For $\ell \in \N$, let $U_{\ell}:J_{k,m} \to J_{k,m\ell^2}$ be the index raising operator from~\cite{ez85}.  
The simplest to define is~$\eta_{\ell}$.  
For $f \in S_k\left( K(N)  \right)^{\epsilon}$ with $\FJ(f) =\sum_{m=1}^{\infty} \phi_m\, \xi^{Nm}$, define 
$A(f)= \sum_{m=1}^{\infty} (\phi_m\vert U_{\ell}) \xi^{N\ell^2 m} \in \FMcusp(k,N\ell^2)$.  
One directly checks $A(f) \in \FMcusp(k,N\ell^2, \epsilon)$ so that by the main result there is an 
$\eta_{\ell}f \in S_k\left( K(N\ell^2)  \right)^{\epsilon}$ with $\FJ(\eta_{\ell}f )= A(f)$.  
This defines $\eta_{\ell}: S_k\left( K(N)  \right) \to  S_k\left( K(N\ell^2)  \right)$ since the entire space 
is the direct sum of the plus and minus forms.  

In order to define $\theta_{\ell}$ and   $\theta'_{\ell}$, set 
$\Xi_{\ell} f=\sum_{m=1}^{\infty} (\phi_m\vert V_{\ell}) \xi^{N\ell m} \in \FMcusp(k,N\ell)$, 
and 
$\Xi'_{\ell} f=\sum_{m=1}^{\infty} \left( \sum_{\delta \mid \gcd(\ell,m)} \delta^{k-1} \phi_{\frac{\ell m}{\delta^2}} \vert U_{\delta} \right) \xi^{N\ell m} \in \FMcusp(k,N\ell)$.  
One directly checks 
$B^{\pm}(f)= \Xi_{\ell}( f) \pm \Xi'_{\ell}( f\vert \mu_N) \in \FMcusp(k,N\ell, \pm 1)$ so that by the main result 
there are $b^{\pm}(f) \in S_k\left( K(N\ell)  \right)^{\pm}$ with 
$\FJ\left( b^{\pm}(f) \right) = B^{\pm}(f)$.  
Defining $\theta_{\ell}f = \frac12(b^{+}(f)+b^{-}(f)) \in S_k\left( K(N\ell)  \right)$ 
and correspondingly 
$\theta'_{\ell}f = \frac12(b^{+}(f\vert \mu_N)-b^{-}(f \vert \mu_N)) \in S_k\left( K(N\ell)  \right)$, 
we have  
$\theta_{\ell},  \theta'_{\ell}: S_k\left( K(N)  \right) \to  S_k\left( K(N\ell)  \right)$ 
with $\FJ\left( \theta_{\ell}f \right) = \Xi_{\ell} f$, 
$\FJ\left( \theta'_{\ell}f \right) = \Xi'_{\ell} f$, and 
$(\theta_{\ell}f)\vert \mu_{N\ell}= \theta'_{\ell}(f \vert \mu_N)$.  
Since these constructions work for  $\ell \in \N$, 
this gives a generalization of the global level raising operators.  
\smallskip

The theory of formal series for arithmetic subgroups developed in~\cite{bruinierraum24} 
frames different hypotheses than we do in the case of the paramodular groups.  
Still, in their Theorem~{4.8}, Bruinier and Raum considered  the implications of their general 
theory for paramodular groups when a single formal series at the standard $1$-cusp is given.  
In our notation, they proved the following.  
For squarefree~$N$, let $\Gamma_0(N)^{*}$ be the extension of~$\Gamma_0(N)$ by all the 
Atkin-Lehner involutions.  
Let $\chi: \Gamma_0(N)^{*} \to \{\pm 1\}$ be a character trivial on~$\Gamma_0(N)$.  
If a formal series $\fracf \in \FM(k,N)$ satisfies 
$\jj(\sigma^*) \AFS(\fracf)= \chi(\sigma) \AFS(\fracf)$ for all $\sigma \in  \Gamma_0(N)^{*}$, 
then~$\fracf$ converges on~$\Half_2$.   
Using simple generators for $\Gamma_0(N)^{*}$, this reproves the cases~$N=2,3$ from~\cite{ipy13}.  
In comparison, our main result allows general~$N$ and assumes a single symmetry under the 
paramodular Fricke involution.
\bigskip

Theorem~\ref{theoremord} only needed to be proven for cuspidal quotients~$\fracf$ 
because the Main Theorem~\ref{thmain3} was reduced to the case of formal series of 
Jacobi cusp forms.  With the Main Theorem in hand, however, we may improve Theorem~\ref{theoremord} 
to obtain a general criterion for the divisibility of paramodular Fricke eigenforms that is interesting 
in its own right because it avoids all discussion of divisors.  

\begin{theorem}
\label{theoremordconclusion}
Let $k,k_0,N \in \N$ and $\epsilon, \epsilon_0 \in \{ \pm 1\}$.  
Let $\fracf \in \FM(k,N)$, $f_0 \in M_{k_0}\left( \KN \right)^{\epsilon_0}$ be nontrivial,  and 
$G \in M_{k_0+k}\left( \KN \right)^{\epsilon_0\epsilon}$.  
If we have $\FJ(G)=  \FJ(f_0) \fracf$ in $\FM(k_0+k,N)$ then 
the meromorphic~$G/f_0$ is holomorphic and $G/f_0 \in M_{k}\left( \KN \right)^{\epsilon}$.   
\end{theorem}
\begin{proof}
We have $\fracf \in \FM(k,N, \epsilon)$ if and only if 
$\jj(F_N^*) \AFS(\fracf) = \epsilon \AFS(\fracf)$.  
Applying $\AFS$ to $\FJ(G)=  \FJ(f_0) \fracf$ we obtain 
$\FS(G)=  \FS(f_0) \AFS(\fracf)$.  
Applying $\jj(F_N^*)$ we have $\epsilon_0 \epsilon\,\FS(G)=  \epsilon_0\FS(f_0)\, \jj(F_N^*)\AFS(\fracf)$.  
Since $\ffs$ is an integral domain and the element $\FS(f_0)$ is nontrivial, 
 we have 
$\jj(F_N^*)\AFS(\fracf) = \epsilon \AFS(\fracf)$, and thus $\fracf \in \FM(k,N, \epsilon)$.  
By the Main Theorem~\ref{thmain3} there is an $f \in M_k\left( \KN \right)^{\epsilon}$ 
such that $\FJ(f)=\fracf$.  
Thus the quotient $G/f_0 =f\in M_{k}\left( \KN \right)^{\epsilon}$ is holomorphic.   
\end{proof}

\bibliographystyle{plain}
\bibliography{bibola}

\end{document}